\documentclass[12pt,a4paper]{amsart}
\usepackage{amsmath}
\usepackage{mathabx}
\usepackage{amstext}
\usepackage{amssymb,mathtools}
\usepackage{array}
\usepackage{amsthm}
\usepackage{dsfont}
\usepackage[latin1]{inputenc}
\usepackage[T1]{fontenc}
\usepackage{mathtools}
\usepackage{appendix}
\usepackage{xcolor}
\usepackage{ulem}
\usepackage[arrow, matrix, curve]{xy}
\usepackage{verbatim}
\usepackage{graphicx}
\usepackage{subfigure}
\usepackage{hyperref}
\usepackage{cases}
\usepackage{esvect}
\usepackage[update,prepend]{epstopdf}

\newtheorem{theorem}{Theorem}
\theoremstyle{plain}
\newtheorem{definition}[theorem]{Definition}
\newtheorem{assumption}[theorem]{Annahme}
\newtheorem{lemma}{Lemma}

\newtheorem{proposition}[theorem]{Proposition}
\theoremstyle{definition}
\newtheorem{example}{Example}

\numberwithin{equation}{section}

\def\eps{\varepsilon}

\def\ri{{\rm i}}

\def\util{\tilde{u}}

\def\Atil{\tilde{A}}
\def\Abartil{\tilde{\bar{A}}}
\def\Ahat{\widehat{A}}
\def\Ahatbar{\widehat{\widebar{A}}}
\def\Uvec{\vec{U}}
\def\Uvecext{\vec{U}^\text{ext}}
\def\Uvecapp{\vec{U}^\text{app}}

\def\uapp{U_\text{app}}

\def\Evec{\vec{E}}

\def\hot{\text{h.o.t.}}
\def\LCME{L_{\text{CME}}}

\def\bi{\begin{itemize}}
\def\ei{\end{itemize}}

\newcommand{\R}{\mathbb{R}}
\newcommand{\C}{\mathbb{C}}
\newcommand{\B}{\mathbb{B}}

\newcommand{\N}{\mathbb{N}}

\newcommand{\Z}{\mathbb{Z}}

\newcommand{\T}{\mathbb{T}}

\newcommand{\cD}{{\mathcal D}}

\newcommand{\cL}{{\mathcal L}}

\newcommand{\cT}{{\mathcal T}}

\newcommand{\cX}{{\mathcal X}}

\renewcommand{\phi}{\varphi}

\newcommand{\pa}{\partial}

 \def\dd{\, {\rm d}}

\DeclareOldFontCommand{\it}{\normalfont\itshape}{\mathit}

\newcommand{\bspm}{\left(\begin{smallmatrix}}\newcommand{\espm}{\end{smallmatrix}\right)}
\newcommand{\bpm}{\begin{pmatrix}}\newcommand{\epm}{\end{pmatrix}}

\def\bs{\begin{satz}}\def\es{\end{satz}}
\def\blem{\begin{lemma}}\def\elem{\end{lemma}}
\def\bthm{\begin{theorem}}\def\ethm{\end{theorem}}
\def\bcor{\begin{corollary}}\def\ecor{\end{corollary}}
\def\beq{\begin{equation}}\def\eeq{\end{equation}}
\def\beqq{\begin{equation*}}\def\eeqq{\end{equation*}}
\def\bal{\begin{align}}\def\eal{\end{align}}
\def\bpf{\begin{proof}}\def\epf{\end{proof}}
\def\bex{\begin{example}}\def\eex{\end{example}}
\def\brem{\begin{remark}}\def\erem{\end{remark}}
\def\bass{\begin{assumption}}\def\eass{\end{assumption}}
\def\bprop{\begin{proposition}}\def\eprop{\end{proposition}}
\def\bdefi{\begin{definition}}\def\edefi{\end{definition}}

\DeclareFontEncoding{FMS}{}{}
\DeclareFontSubstitution{FMS}{futm}{m}{n}
\DeclareFontEncoding{FMX}{}{}
\DeclareFontSubstitution{FMX}{futm}{m}{n}
\DeclareSymbolFont{fouriersymbols}{FMS}{futm}{m}{n}
\DeclareSymbolFont{fourierlargesymbols}{FMX}{futm}{m}{n}
\DeclareMathDelimiter{\VERT}{\mathord}{fouriersymbols}{152}{fourierlargesymbols}{147}

\def\bi{\begin{itemize}}
\def\ei{\end{itemize}}
\def\ben{\begin{enumerate}}
\def\een{\end{enumerate}}
\usepackage{enumitem}
\setitemize{leftmargin=*}

\begin{document}

\title[Coupled Mode Equations and Gap Solitons in Higher Dimensions]{Coupled Mode Equations and Gap Solitons in Higher Dimensions}
\author{Tom\'a\v{s} Dohnal$^{1}$ and Lisa Wahlers$^{2}$}
\address{$^{1}$ Martin-Luther-Universit\"at Halle-Wittenberg, Institut f\"ur Mathematik, D-06099 Halle (Saale), Germany; $^{2}$ Technische Universit\"at Dortmund, Fakult\"at
      f\"ur Mathematik, Vogelpothsweg 87, D-44227 Dortmund, Germany}
\email{tomas.dohnal@mathematik.uni-halle.de, lisa.wahlers@math.tu-dortmund.de}
\date{\today}

%\makeindex

\begin{abstract}
We study waves-packets in nonlinear periodic media in arbitrary ($d$) spatial dimension, modeled by the cubic Gross-Pitaevskii equation. In the asymptotic setting of small and broad waves-packets with $N\in \N$ carrier Bloch waves the effective equations for the envelopes are first order coupled mode equations (CMEs). We provide a rigorous justification of the effective equations. The estimate of the asymptotic error is carried out in an $L^1$-norm in the Bloch variables. This translates to a supremum norm estimate in the physical variables. In order to investigate the existence of gap solitons of the $d$-dimensional CMEs, we discuss spectral gaps of the CMEs. For $N=4$ and $d=2$ a family of time harmonic gap solitons is constructed formally asymptotically and numerically. Moving gap solitons have not  been found for $d>1$ and for the considered values of $N$ due to the absence of a spectral gap in the standard moving frame variables. 
\end{abstract}

\maketitle
%\thispagestyle{empty}
%\newpage
%\cleardoublepage
%\thispagestyle{empty}
%\tableofcontents
%\setcounter{chapter}{-1}
%\setcounter{secnumdepth}{3}
%\setlength\parindent{0pt}

   \vspace*{2mm} {\bf Key-words:} waves-packets, coupled mode equations, periodic media, gap soliton, approximation error, Gross-Pitaevskii 

     \vspace*{2mm}
		
\setcounter{page}{1}

%Introduction
\section{Introduction}
\label{S:intro}

One of the intriguing properties of nonlinear dispersive equations is the possibility of coherent spatially localized waves. We study such waves in the form of waves-packets in periodic media of arbitrary dimension. There are two standard asymptotic scalings of such waves-packets: waves-packets with a single carrier Bloch wave and those with multiple ($N$) carrier Bloch waves with different group velocities. In the former case the waves-packet propagates at a velocity asymptotically close to the carrier group velocity ($v_g$) and the waves-packet envelope can be (in generic situations) described by the constant coefficient nonlinear Schr\"odinger equation \cite{BSTU06,GWH01,DR18}. This is based on Taylor-expanding the band structure at the carrier wavevector to second order  and using the moving frame variable with the velocity $v_g$. In the latter case $N$ carriers with the same temporal frequency $\omega_0$ are used and the band structure is expanded only to first order producing a system of $N$ first order amplitude equations for the $N$ envelopes, so called coupled mode equations (CMEs) \cite{AW89,SU01,GWH01,AP05,DH17,Pelinov_2011}. In the one dimensional case the CMEs have a family of solitary waves parametrized by the velocity $v\in (-1,1)$ (after a rescaling). As these solitary waves exist in the spectral gap of the CMEs, they are called gap solitons. In a special case (corresponding to an infinitesimal contrast in the periodicity in \eqref{eq:PNLS-W}) there is even an explicit gap soliton family. In the general case only numerical representations of gap solitons have been found  \cite{D14,DH17}. An important question is whether also in the higher dimensional ($d$-dimensional) CMEs such a family of gap solitons exists. These gap solitons would then approximate coherent waves-packets of the original model with a range of velocities in a $d$-dimensional interval and with frequencies asymptotically close to $\omega_0$. Such a tunability is interesting from a phenomenological as well as an applied point of view.

As a prototypical model we consider in this paper the $d-$dimensional periodic Gross-Pitaevskii equation with the cubic nonlinearity
\beq\label{eq:PNLS-W}
\ri \pa_t u +\Delta u -(V(x)+\eps W(x)) u -\sigma(x) |u|^2u=0, \ x\in \R^d,t\in \R,
\eeq
where 
$$V,\sigma \in C(\R^d,\R) \text{ are } 2\pi\Z^d-\text{periodic},$$
$\eps >0$ is a small parameter, and 
$$W(x)=\sum_{m=-m_*}^{m_*}a_m e^{\ri l^{(m)}\cdot x}$$ 
is a real periodic function specified in (A4) below. Precise regularity conditions on $V$ and $\sigma$ appear in (A5)  and are needed for the justification of amplitude equations in Sec. \ref{S:justif}. Our analysis can be carried over to other cubically nonlinear wave models with only small modifications. An example is the wave equation $\pa_t^2 u -\Delta u +(V(x)+\eps W(x))u+\sigma(x)u^3=0.$

Our formal asymptotic ansatz for an $N-$mode waves-packet solution to \eqref{eq:PNLS-W} is
\begin{alignat}{2}
u^{\text{app}}(x,t)\coloneqq \varepsilon^{1/2} \sum\limits_{j=1}^N A_j(\varepsilon x, \varepsilon t)p_j(x) e^{\ri (k^{(j)} \cdot x-\omega_0 t)} , \label{eq:uapp}
\end{alignat} 
where $p_j(x) e^{\ri (k^{(j)} \cdot x-\omega_0 t)}, j=1,\dots,N$ are the carrier Bloch waves, see Sec. \ref{S:bloch-form-deriv}. The functions $p_j(x):=p_{n_j}(x,k^{(j)})$ are eigenfunctions of the Bloch eigenvalue problem \eqref{E:Bl-ev} and the wave-vectors $k^{(j)}$ lie in the Brillouin zone 
$$\B:=(-\tfrac{1}{2},\tfrac{1}{2}]^d$$ 
corresponding to the $2\pi\Z^d$ periodicity. The asymptotic parameter $\eps$ in the ansatz is the same as the $\eps$ in equation \eqref{eq:PNLS-W}. The role of the periodic perturbation $\eps W$ of the linear potential is to provide for linear coupling among the modes. This is explained below and in more detail in Sec. \ref{S:bloch-form-deriv}. Note that if the periods of $V$ and $W$ are not commensurate, the potential $V+\eps W$ is quasiperiodic. In our asymptotic result this aspect plays, however, no role. 

Under assumptions (A1), (A2), and (A4) one can formally derive the following coupled mode equations (CMEs) for the envelopes $A_j,j=1,\dots,N$. For details see Sec. \ref{S:bloch-form-deriv}.
\beq
\begin{aligned}
\ri (\pa_TA_j+v_g^{(j)}\cdot \nabla_X A_j) + \sum_{r=1}^N\kappa_{jr} A_r + N_j(\vec{A})=0, \ j = 1,\dots,N,
\label{E:CME}
\end{aligned}
\eeq
where for $j,r\in \{1,\dots,N\}$
\begin{align*}
N_j(\vec{A})&:=\sum_{\stackrel{(\alpha,\beta,\gamma)\in\{1,\dots,N\}^3}{k^{(\alpha)}-k^{(\beta)}+k^{(\gamma)}\in k^{(j)}+\Z^d}} \gamma_j^{(\alpha,\beta,\gamma)}A_{\alpha}\overline{A}_{\beta}A_{\gamma},\\ 
\gamma_j^{(\alpha,\beta,\gamma)} &:= -\int_{\T} \sigma(x) p_{\alpha}(x)\overline{p_{\beta}}(x)p_{\gamma}(x)\overline{p_j}(x)e^{\ri (k^{(\alpha)}-k^{(\beta)}+k^{(\gamma)}-k^{(j)})\cdot x} \dd x,\\
\kappa_{jr}&:=-\sum_{\stackrel{m\in\{-m_*,\dots,m_*\}}{k^{(r)}+l^{(m)}\in k^{(j)}+\Z^d}}a_{m}\int_\T e^{\ri(k^{(r)}+l^{(m)}-k^{(j)})\cdot x} p_r(x)\overline{p_j}(x)\dd x,
\end{align*}
and where  $\T:=\R^d/(2\pi \Z^d)$ is a $d$-dimensional torus. Due to the assumptions on $a_m$ and $l^{(m)}$ in (A4) the matrix $\kappa=(\kappa_{jr})_{j,r=1}^N$ is Hermitian.

In the three dimensional case CMEs have been formally derived in \cite{AP05} as envelope equations for the Maxwell equations. In \cite{GMS08} the $d-$dimensional CMEs \eqref{E:CME} have been analyzed. However, because the authors start with \eqref{eq:PNLS-W} with $W\equiv 0$, the coupling coefficient matrix $\kappa$ vanishes. With $\kappa = 0$ there is no spectral gap of the linear part of CMEs because the dispersion relation of the CMEs is then given by the four hyperplanes $K\mapsto v_g^{(j)}\cdot K, K \in \R^d$.
No (exponentially) localized solitary waves of CMEs are thus expected with $\kappa=0$. We show in Sec. \ref{S:CMEs_band-str} that at least for $N=4$ there are $\kappa\neq 0$ such that a spectral gap exists and standing solitary waves can be constructed.

One of the main contributions of the paper is a rigorous justification of \eqref{E:CME} as valid amplitude equations for the waves-packets \eqref{eq:uapp} in \eqref{eq:PNLS-W}. In \cite{GMS08} such a justification was performed for $W\equiv 0$. The authors work in the space $H^s_\eps(\R^d):=\{f\in L^2(\R^d):\int_{\R^d}(1+|\eps k|^2)|\hat{f}(k)|^2\dd k<\infty\}$. Because of the loss of $d/2$ powers of $\eps$ when evaluating the $L^2(\R^d)$ norm of a function $f(\eps \cdot)$, the authors are forced to work with additional higher order terms in the asymptotic expansion of the waves-packet. In order to invert the linear operator near the concentration $k$-points of the new correction terms the authors of \cite{GMS08} impose a closed mode system condition (a condition on the terms generated by a repeated application of the nonlinearity on the ansatz). In our approach we work in an $L^1$ space in the Bloch variables similarly to \cite{DH17}, where the one dimensional case was considered. The $L^1$-approach for the justifictaion of coupled mode equations was used already in \cite{SU01} and \cite{Pelinov_2011}, where the case of periodicity with an infinitesimal contrast (i.e. $V\equiv 0$ in \eqref{eq:PNLS-W}) was considered and the analysis was carried out in Fourier variables. The one dimensional CMEs for the infinitesimal contrast were justified also in \cite{GWH01} with estimates in Sobolev norms. In \cite{DH17} the one dimensional case with $V\neq 0$ is treated and Bloch variables are needed -similarly to the current paper. By an analog of the Riemann-Lebesgue Lemma we can estimate the supremum of the error in $x$-variables in terms of the $L^1$ norm in Bloch variables. Via this approach we avoid working with high order terms in the asymptotics.

We use the following assumptions.
\bi
\item[(A1)] $k^{(1)},\dots, k^{(N)}\in \B$ and $n_1,\dots,n_N\in \N$ are such that $(n_1,k^{(1)}),\dots,(n_N,k^{(N)})$ are pairwise distinct and $\omega_{n_j}(k^{(j)}), j=1,\dots,N$ are eigenvalues of of \eqref{E:Bl-ev} with 
$$\omega_{n_1}(k^{(1)})=\omega_{n_2}(k^{(2)})=\dots=\omega_{n_N}(k^{(N)})=\omega_0\in \R.$$
If $k^{(j_1)}=k^{(j_2)}=\dots=k^{(j_I)}$ for some $\{j_1,\dots,j_I\}\subset \{1,\dots,N\}$, then the corresponding eigenfunctions $p_{n_{j_1}}(x,k^{(j_1)}),\dots,$ $p_{n_{j_I}}(x,k^{(j_I)})$ are $L^2(\T)$-orthogonal.
\item[(A2)] $v_g^{(j)}$ is twice continuously differentiable at $k^{(j)}$ for each $j=1,\dots,N$. 
\item[(A3)] $|\omega_{n}(l)-\omega_0|>\delta>0$ for all $(n,l)\in (\N\times J) \setminus \{(n_1,k^{(1)}), \dots,(n_N,k^{(N)})\}$ with $J$ defined in \eqref{E:J}.
\item[(A4)] $W(x)=\sum_{m=-m_*}^{m_*} a_m e^{\ri l^{(m)}\cdot x}$, where $m_*\in \N$, the vectors $l^{(m)}\in \R^d,m =-m_*,\dots,m_*$ are pairwise distinct, and 
$$ a_{-m}=\overline{a_m}, l^{(-m)}=-l^{(m)} \text{ for all } m\in\{0,1,\dots,m_*\}.$$ 
\item[(A5)] $V \in H^{2\left\lceil \tfrac{d}{2}\right\rceil +d+\delta} \left( \T, \R \right), \sigma \in H^{2\left\lceil \tfrac{d}{2}\right\rceil +2}\left( \T, \R \right), \delta>0$.
\item[(A6)] For each $j\in\{1,\dots,N\}$ there is $L>0$, such that
\beq\label{E:Lipschitzpn}
\|p_{n_j}(\cdot,k)-p_{n_j}(\cdot,k^{(j)})\|_{L^2(\T)}\leq L |k-k^{(j)}|
\eeq
for all $k$ in a neighborhood of $k^{(j)}$. 
\ei
For $a\in \R$ we denote by $\lceil a \rceil$ the smallest integer larger than or equal to $a$.

Note that if at each $k^{(j)}$ the eigenvalue $\omega_{n_j}(k^{(j)})$ is simple, then \eqref{E:Lipschitzpn} holds automatically, see \cite[Sec. VII.2]{Kato_1995}.

The potential $\eps W$ is introduced in order to linearly couple the modes of the ansatz \eqref{eq:uapp}. In detail, two modes $e^{\ri k^{(j)}\cdot x}p_{n_j}(x,k^{(j)})$ and $e^{\ri k^{(\tilde{j})}\cdot x}p_{n_{\tilde{j}}}(x,k^{(\tilde{j})})$ with $j,\tilde{j} \in \{1,\dots,N\}$ are coupled if 
$$k^{(j)}-k^{(\tilde{j})}=l^{(m)}+\mu \text{ for some } m\in\{-m_*,\dots,m_*\}, \mu\in \Z^d,$$
and
$$a_{m}\int_\T e^{\ri(k^{(\tilde{j})}+l^{(m)}-k^{(j)})\cdot x} p_{\tilde{j}}(x)\overline{p_j}(x)\dd x \neq 0.$$
The first condition implies that  
$$e^{\ri k^{(\tilde{j})}\cdot x}e^{\ri l^{(m)}\cdot x} = e^{\ri k^{(j)}\cdot x}\psi(x),$$
where $\psi$ is $2\pi\Z^d$-periodic. For details see Section \ref{S:bloch-form-deriv}.

%internal note: the above conditions imply a_0\in \R (necessary for W real); The presence of the zero mode in W is needed for the linear coupling if k^{(j_1)}=k^{(j_2)} for some j_1\neq j_2

Our justification result is 
\bthm \label{T:justif}
Assume (A1)-(A6) and let $\vec{A}$ be a solution of \eqref{E:CME} with $\hat{A}_j \in C\left( [0,T_0], L^1_{s_A}(\R^d)\right.$ $\left.\cap  L^2 (\R^d)\right)$, $\pa_T\hat{A}_j \in C\left( [0,T_0], L^1(\R^d)\right)$ for some $T_0>0$, all $j=1\dots,N$ and some $s_A > 2\lceil \tfrac{d}{2}\rceil +d+2$. Then there are constants $c>0$ and $\eps_0>0$, such that if $u(\cdot,0)=u^{\text{app}}(\cdot,0)$ with $u^\text{app}$ given by \eqref{eq:uapp}, then for all $\eps \in (0, \eps_0)$ the solution $u$ of \eqref{eq:PNLS-W} satisfies $u(x,t) \rightarrow 0$ as $|x| \rightarrow \infty$ and
\begin{alignat}{2}\label{E:est-main}
\|u(\cdot, t)-u^{\text{app}}(\cdot,t)\|_{C_b^0} \leq c \eps^{3/2} \text{ for all } t \in [0,\eps^{-1}T_0].
\end{alignat} 
\ethm

The rest of the paper is organized as follows. In Sec. \ref{S:bloch-form-deriv} we review the Bloch eigenvalue problem and provide a brief formal derivation of the CMEs \eqref{E:CME}.
In Sec. \ref{S:CME-asymp} an asymptotic approximation of  standing (time harmonic) gap solitons of the CMEs is constructed for frequencies in a spectral gap  and asymptotically close to the spectrum. Sec. \ref{S:CMEs_band-str} studies the dispersion relation and the gap structure of CMEs with $N$ odd and $N=2,4$ and finds that for $d\geq 2$ a necessary condition for the existence of a standing gap soliton is $N\in \{4,6,8, \dots\}$. An example with a gap for $N=4$ is produced and in Sec. \ref{S:num-stand} a standing gap soliton is found numerically. The question of moving gap solitons is addressed in Sec. \ref{S:mov_solit}. It is shown that for $d\geq 2$ and $N\leq 4$ no spectral gap is present in the CMEs in moving frame variables. Finally, Sec. \ref{S:justif} presents a proof of Theorem \ref{T:justif}, i.e. a rigorous justification of CMEs \eqref{E:CME} as amplitude equations for \eqref{eq:PNLS-W}.
%-----------------------------------------------------------------
\subsection{The Band Structure and a Formal Derivation of the Coupled Mode Equations}\label{S:bloch-form-deriv}

The bounded elementary solutions of the linear part of \eqref{eq:PNLS-W} at $\eps=0$  are the Bloch waves
\beq\label{E:Bloch-wv}
p_n(x,k)e^{\ri (k\cdot x-\omega_n(k) t)}, k \in \B, n \in \N,
\eeq
where $(\omega_n(k),p_n(x,k)), n \in \N$ are eigenpairs of the Bloch eigenvalue problem
\beq\label{E:Bl-ev}
\begin{aligned}
&\mathcal{L}(x,k) p_n(x,k) = \omega_n(k) p_n(x,k), \ x\in \T,\\
&\mathcal{L}(x,k) := -|\nabla + \ri k|^2 + V(x)=-(\Delta  + 2\ri k\cdot \nabla - |k^2|)+ V(x).
\end{aligned}
\eeq
The graph $(k,\omega_n(k))_{n\in \N}$ is called the band structure. The eigenfunctions $p_n(\cdot,k)$ are automatically $2\pi$-periodic in each coordinate. Because of the self-adjointness of $\mathcal{L}(\cdot,k)$ the eigenfunction sequence $(p_n(\cdot,k))_{n\in \N}$ can be chosen $L^2(\T)$-orthonormal. 

We have the following asymptotics of the eigenvalues. There exist constants $c_1,c_2>0$ such that for all $n\in \N, k \in \B$
\beq\label{E:as-distr}
c_1 n^{2/d} \leq |\omega_n(k)| \leq c_2 n^{2/d},
\eeq
, see \cite[p.55]{Hoerm_85}.

The group velocity of the Bloch wave \eqref{E:Bloch-wv} is $\nabla \omega_n(k)$, which has got also the following integral representation obtained by differentiating \eqref{E:Bl-ev} in $k$
$$\nabla \omega_n(k)=-2\langle (\nabla + \ri k)p_n(\cdot,k), p_n(\cdot,k)\rangle_{L^2(\T)}.$$

For the operator $-\Delta +V:L^2(\R^d)\to L^2(\R^d)$ it is $\text{spec}(-\Delta+V)=\cup_{n\in \N}\omega_n(\B)$, see 
Chapter 3 in \cite{DLPSW_2011}. 

\begin{figure}[ht!]
\includegraphics[scale=0.6]{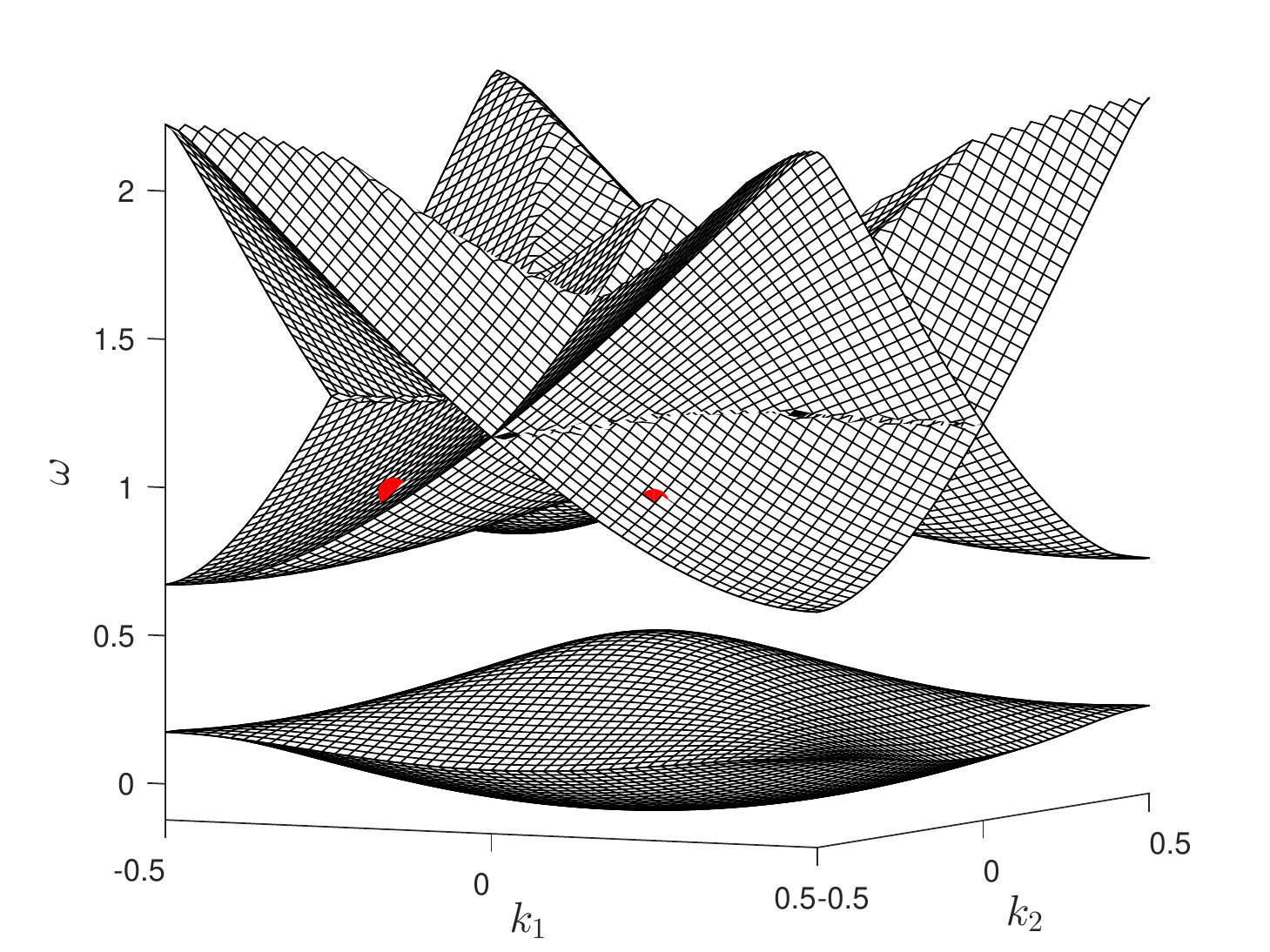}
\caption{The band structure of \eqref{E:Bloch-wv} for $V(x)=\cos(x_1)\cos(x_2)$. The points $(k^{(1)},\omega_0)$ and $(k^{(2)},\omega_0)$ are marked by red dots. The points $k^{(3)}$ and $k^{(4)}$ satisfy $k^{(3)}=-k^{(1)}$ and $k^{(4)}=-k^{(2)}$.}
\label{F:band-str}
\end{figure}

\begin{figure}[ht!]
\includegraphics[scale=0.6]{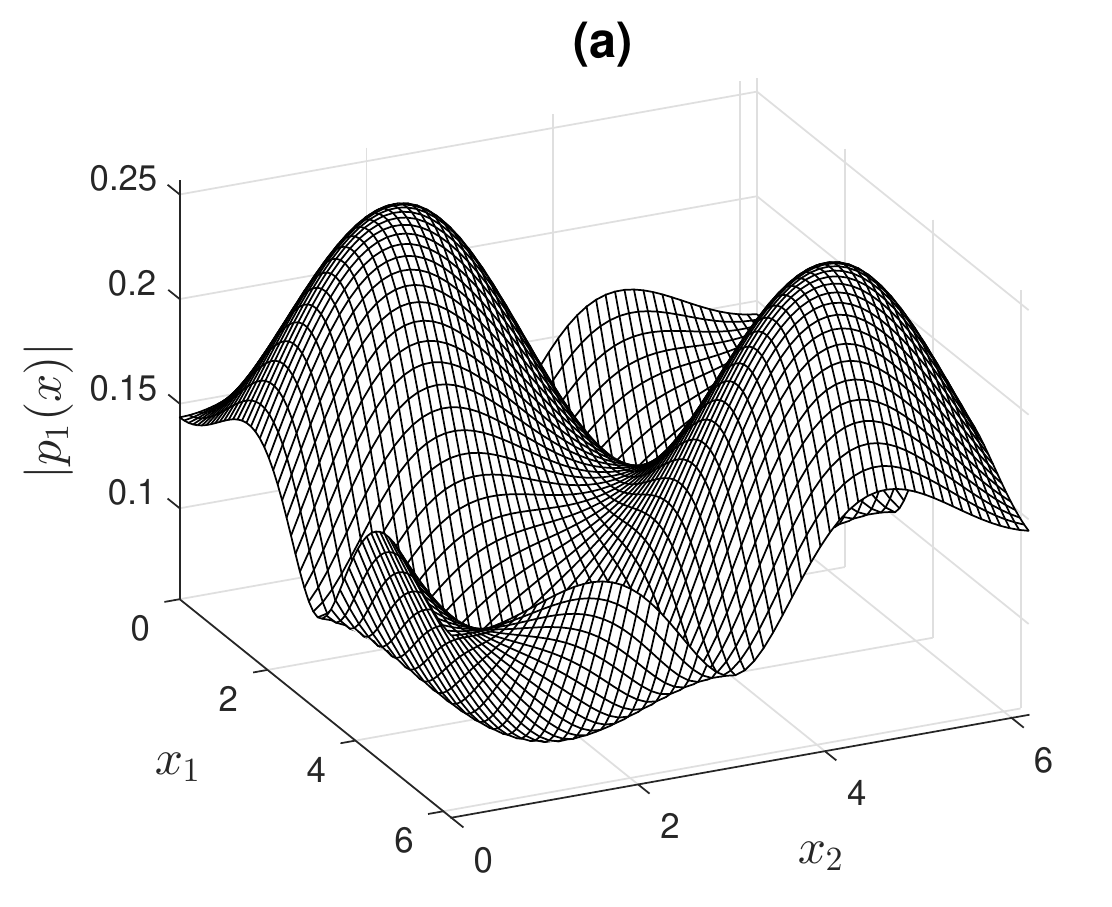}
\includegraphics[scale=0.6]{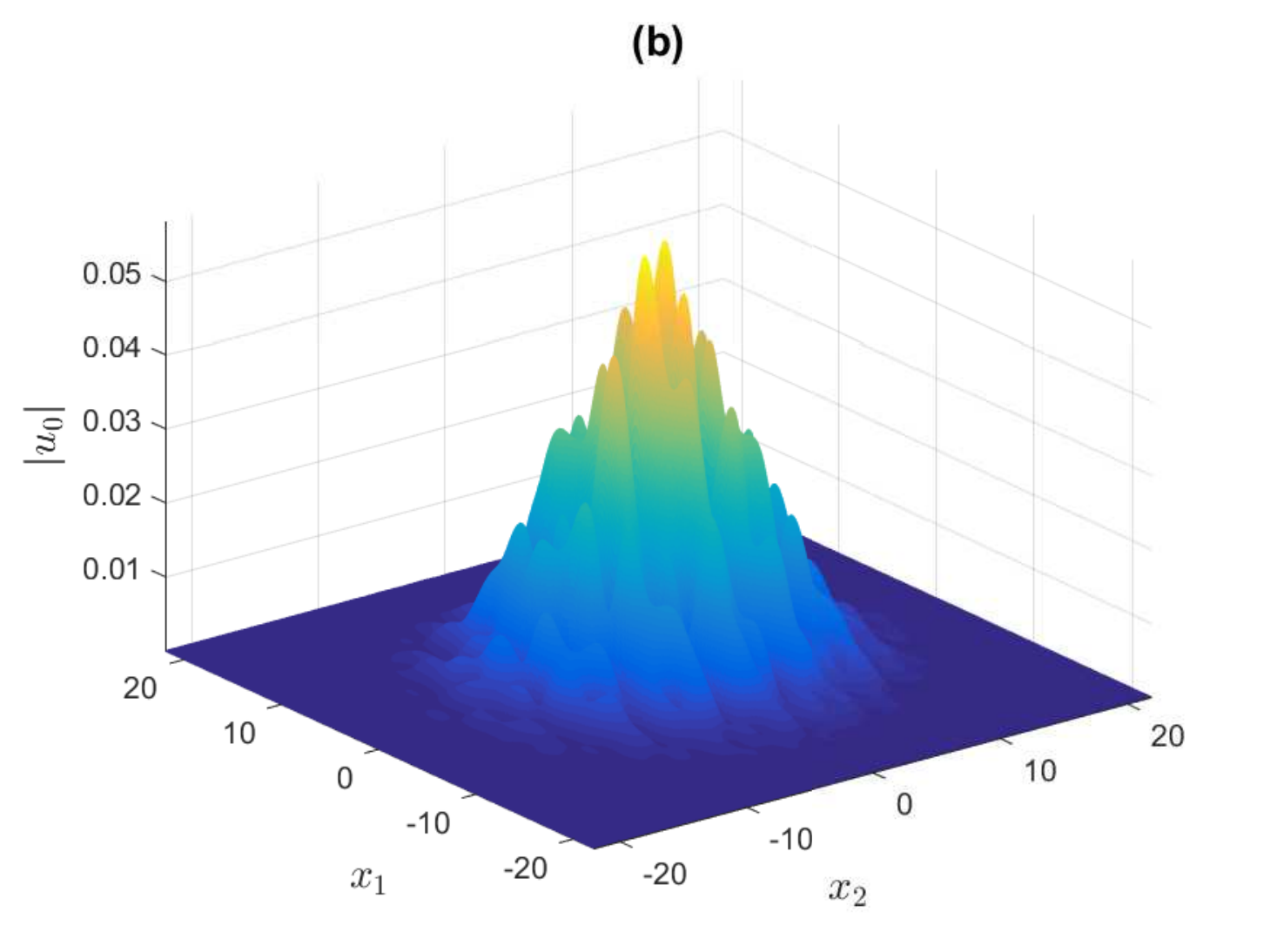}
\caption{(a) The modulus of $p_1:=p_{n_1}(\cdot,k^{(1)})$, where $n_1=4, k^{(1)}=(-0.2,-0.4)$, see Fig. \eqref{F:band-str}. (b) The modulus of the approximate ansatz $u^\text{app}(\cdot,0)$ for $N=4$, the points $k^{(j)}, j=1,\dots,4$ as in Fig. \eqref{F:band-str}, $\eps=0.12$, and $A_1(X,0)=\dots=A_4(X,0)=e^{-|X|^2}$.}
\label{F:p1_u0}
\end{figure}

Figure \ref{F:band-str} shows the  band structure $k\mapsto (\omega_n(k))_{n\in \N}$ (with $k\in \B$) for the example $V(x)=\cos(x_1)\cos(x_2)$. Only the first six eigenvalues $k\mapsto \omega_n(k), n=1,\dots,6$ are plotted. To illustrate the asymptotic ansatz \eqref{eq:uapp}, we choose the example $N=4$ with $k^{(1)}=(-0.2,-0.4),k^{(2)}=(0.2,-0.4)$, $k^{(3)}=-k^{(1)}, k^{(4)}=-k^{(2)}$ and $n_1=\dots=n_4=4$. With the potential $V(x)=\cos(x_1)\cos(x_2)$ we have $\omega_0:=\omega_4(k^{(1)})=\dots = \omega_4(k^{(4)})\approx 0.9942$. The points $k^{(1)}$ and $k^{(2)}$ are marked in Fig. \ref{F:band-str}. In Fig. \ref{F:p1_u0} (a) we plot $|p_1|$, where $p_j:=p_{n_j}(\cdot,k^{(j)}), j=1,\dots,4$. The other Bloch functions are related by $p_2(x)=p_1(-x_1,x_2), p_3=\overline{p_1}$ and $p_4=\overline{p_2}$. Finally, Fig. \ref{F:p1_u0} (b) shows $|u^\text{app}(\cdot,0)|$ for the above setting and with $A_1(X,0)=\dots=A_4(X,0)=e^{-|X|^2}, \eps=0.12$. 

Ansatz \eqref{eq:uapp} for an approximate waves-packet uses $N$ elementary Bloch waves as carriers. Due to the slow dependence on $x$ of the envelope $A_j(\eps x,\eps t)$ the ansatz is concentrated in the Bloch variable $k$ near the points $k=k^{(1)}, \dots, k^{(N)}$, see Sec. \ref{S:justif} for the definition of Bloch variables. In the formal derivation of the CMEs \eqref{E:CME} we plug \eqref{eq:uapp} in \eqref{eq:PNLS-W} and collect terms of (formally) equal orders of $\eps$. All $O(\eps^{1/2})$ terms vanish due to the fact that $(\omega_0,p_j)$ are eigenpairs of \eqref{E:Bl-ev} at $k=k^{(j)}, j=1,\dots,N$. At $O(\eps^{3/2})$ both the nonlinearity and the term $\eps Wu^\text{app}$ contribute. We consider in the formal part only terms proportional to $e^{\ri k^{(j)}\cdot x}, j=1,\dots,N$ times a $2\pi \Z^d$-periodic function, i.e. only terms concentrated in the Bloch variable near $k^{(j)}, j=1\dots,N$. For each $j \in \{1,\dots,N\}$ we get
\beq\label{E:res-uapp}
\begin{aligned}
&\sum_{\stackrel{r\in \{1,\dots,N\}}{k^{(r)}=k^{(j)}}}\left(\ri(\pa_T A_r+v_g^{(r)}\cdot \nabla A_r)p_r(x)e^{\ri k^{(j)}\cdot x}\right) -\sum_{r=1}^N A_r\sum_{\stackrel{m\in\{-m_*,\dots,m_*\}}{k^{(r)}+l^{(m)}\in k^{(j)}+\Z^d}}a_{m}e^{\ri(k^{(r)}+l^{(m)})\cdot x} p_r(x)\\
&-\sigma(x) \sum_{\stackrel{(\alpha,\beta,\gamma)\in\{1,\dots,N\}^3}{k^{(\alpha)}-k^{(\beta)}+k^{(\gamma)}\in k^{(j)}+\Z^d}} A_{\alpha}\overline{A}_{\beta}A_{\gamma} p_{\alpha}(x)\overline{p_{\beta}}(x)p_{\gamma}(x)e^{\ri (k^{(\alpha)}-k^{(\beta)}+k^{(\gamma)})\cdot x}. 
\end{aligned}
\eeq
A necessary condition for the smallness of the overall residual for $u^\text{app}$ is that the $L^2(\T)$-projection of these terms onto $p_j(\cdot)e^{\ri k^{(j)}\cdot }$ vanishes. Because of the orthogonality of $\{p_r: r=1,\dots,N, k^{(r)}=k^{(j)}\}$ (see (A1)), this condition is exactly equation \eqref{E:CME}.

Of course, ansatz \eqref{eq:uapp} with $\vec{A}$ satisfying \eqref{E:CME} does not yield a zero residual in \eqref{eq:PNLS-W} because we omitted in $Wu^\text{app}$ and in $|u|^2u$ terms not proportional to $e^{\ri k^{(j)}\cdot x}, j=1,\dots,N$ times a $2\pi \Z^d$-periodic function and because only the projection of \eqref{E:res-uapp} onto $p_j(\cdot)e^{\ri k^{(j)}\cdot }$ was forced to vanish. In Sec. \ref{S:justif} we prove that, despite this, on time intervals of size $O(\eps^{-1})$ the residual is $O(\eps^{5/2})$ and the error between $u$ and $u^\text{app}$ is $O(\eps^{3/2})$.

%-------------------------------------------------------------------
%-------------------------------------------------------------------
%-------------------------------------------------------------------

\section{Asymptotic Construction of Gap Solitons of CMEs}
\label{S:CME-asymp}

As explained in the introduction, one of the motivations for studying CMEs is the possible existence of gap solitons. 
When inserted into $u^\text{app}$ gap solitons generate nearly solitary wave-packets of \eqref{eq:PNLS-W}. Standing solitary waves of CMEs are of the form $\vec{A}(X,T)=e^{-\ri \Omega T}\vec{B}(X)$, where $\vec{B}$ is spatially localized and satisfies 
\beq\label{E:CME-stand}
\Omega \vec{B} = L_\text{CME}(\nabla)\vec{B} - \vec{N}(\vec{B}),
\eeq
where
$$L_{\text{CME}}(\nabla):=\begin{pmatrix} 
-\ri v_g^{(1)}\cdot \nabla & & \\
& \ddots & \\
& &  -\ri v_g^{(N)}\cdot \nabla 
\end{pmatrix} - \kappa.$$
For $\Omega\notin \text{spec}(\LCME)$ but asymptotically near a spectral edge an approximation via nonlinear Schr\"odinger (NLS) asymptotics is possible. These asymptotics require that the dispersion relation is locally (near the spectral edge) of a parabolic shape. Next, we provide a formal derivation of the asymptotic expansion in order to use this in a numerical construction of a standing gap soliton $\vec{A}(X,T)=e^{-\ri \Omega T}\vec{B}(X)$.

The spectrum $\text{spec}(\LCME)$ is given by the range of $\Omega:\R^d\to \R$ as a solution of the dispersion relation, i.e. of
$$\Omega(K)\in \text{spec}(\LCME(\ri K)), \ K\in \R^d.$$
Let us denote the $N$ eigenpairs of
$$\Omega(K)\vec{\eta}(K)=\LCME(\ri K)\vec{\eta}(K)$$
by $(\Omega_j(K),\vec{\eta}^{(j)}(K)), j=1,\dots,N$ with $\Omega_j:\R^d\to \R, \vec{\eta}^{(j)}:\R^d\to \R^N, \|\vec{\eta}^{(j)}(K)\|_{l^2(\R^N)}=1$. We make the following assumptions
\bi
\item[(A1$_\text{CME}$)] $\Omega_*:=\Omega_{j_0}(K_0)\in \pa(\text{spec}(\LCME))$ for some $j_0\in\{1,\dots,N\}, K_0\in \R^d$,
\item[(A2$_\text{CME}$)] $\Omega_{j_0}(K)$ is a simple  eigenvalue for all $K$ in a neighborhood of $K_0$,
\item[(A3$_\text{CME}$)] the Hessian $D^2\Omega_{j_0}(K_0)$ is definite,
\ei
The vanishing of $\nabla \Omega_{j_0}(K_0)$ is automatic as the edge is an extremal value of $\Omega_{j_0}$. Note that (A3$_\text{CME}$) implies that the extremum of $\Omega_{j_0}(K)$ at $K=K_0$ is isolated and (A2$_\text{CME}$) impies that $K\mapsto \Omega_{j_0}(K)$ is analytic at $K_0$.

 Our asymptotic ansatz for a gap soliton at $\Omega = \Omega_* + \eps^2 \lambda, 0 <\eps \ll 1$, is a slowly varying envelope $\eps C(\eps X)$ modulating the plane wave $e^{\ri K_0\cdot X}\vec{\eta}^{(j_0)}(K_0)$. We define the ansatz in Fourier variables as:
\beq\label{E:ansatzB}
\Omega = \Omega_* + \eps^2 \lambda, \quad \vec{\hat{B}}(K) \sim \vec{\hat{B}}_\text{ans}(K):=\eps^{1-d} \hat{C}\left(\frac{K-K_0}{\eps}\right)\vec{\eta}^{(j_0)}(K),
\eeq
where the scalar $\lambda =O(1)$ and the envelope $C:\R^d\to\R$ are to be determined. Inserting $\Omega$ and $\vec{B}_\text{ans}$ into \eqref{E:CME-stand}, we obtain the residual
$$\eps^{1-d} \left[\Omega_*+\eps^2\lambda-\Omega_{j_0}(K_0+\eps \kappa)\right]\vec{\eta}^{(j_0)}(K_0+\eps \kappa)\hat{C}(\kappa) + \beta(\kappa)\vec{\eta}^{(j)}(K_0+\eps \kappa) + \widehat{\vec{N}}_\perp(\vec{B}_\text{ans})(K_0+\eps \kappa),$$
where 
$$\kappa:=\frac{K-K_0}{\eps}, \quad  \beta(\kappa):=\vec{\eta}^{(j_0)}(K_0+\eps\kappa)^T\widehat{\vec{N}}(\vec{B}_\text{ans})(K_0+\eps\kappa)$$ 
is the coefficient of the nonlinearity component parallel to $\vec{\eta}^{(j_0)}(K)$, and 
$$\widehat{\vec{N}}_\perp(\vec{B}_\text{ans})(K_0+\eps \kappa):= \widehat{\vec{N}}(\vec{B}_\text{ans})(K_0+\eps \kappa) - \beta(\kappa)\vec{\eta}^{(j_0)}(K_0+\eps \kappa)$$
is the orthogonal part. We expand
$$\Omega_{j_0}(K_0+\eps \kappa) = \Omega_* + \tfrac{1}{2}\eps^2 \kappa^T D^2\Omega_{j_0}(K_0) \kappa + O(\eps^3).$$
The $O(\eps^{1-d})$ part of the residual then vanishes. Note that both $\beta$ and $\widehat{\vec{N}}_\perp(\vec{B}_\text{ans})$ are $O(\eps^{3-d})$. 

In order to make the $O(\eps^{3-d})$ part of the residual component parallel to $\eta^{(j_0)}(K)$ vanish, we set 
$$\lambda \hat{C}(\kappa)-\tfrac{1}{2} \kappa^T D^2\Omega_{j_0}(K_0) \kappa \hat{C}(\kappa)+\beta(\kappa)=0, \kappa\in\R^d. $$
This effective equation for $C$ is not very practical due to the complicated structure of $\beta$. However, if we replace in $\beta$ the vector $\vec{\eta}^{(j_0)}(K_0+\eps\kappa)$ by $\vec{\eta}^{(j_0)}(K_0)$, we formally arrive at $\beta(\kappa)=\vec{\eta}^{(j_0)}(K_0)^T\vec{N}(\vec{\eta}^{(j_0)}(K_0))(\hat{C}*\hat{C}*\hat{\overline{C}})(\kappa)$. The resulting effective equation in physical variables is the generalized NLS
\beq\label{E:NLS-eff}
\lambda C+\tfrac{1}{2} \nabla^T D^2\Omega_{j_0}(K_0) \nabla C+\Gamma |C|^2C=0,
\eeq
where $\Gamma := \vec{\eta}^{(j_0)}(K_0)^T\vec{N}(\vec{\eta}^{(j_0)}(K_0))$. The above replacement can be shown to generate only higher order terms in the residual.

The orthogonal component $\widehat{\vec{N}}_\perp(\vec{B}_\text{ans})$ of the $O(\eps^{3-d})$ part of the residual can be eliminated by adding to the ansatz $\vec{\hat{B}}_\text{ans}$ the correction term $\eps^{3-d}\left(Q(L_\text{CME}-\Omega_*)Q\right)^{-1}\widehat{\vec{N}}_\perp(\vec{B}_\text{ans})$ with $Q:=\chi_{B_{\eps^{1/2}}(K_0)}(K)(I-P), Pv:=(\eta^{(j_0)}(K)^Tv)\eta^{(j_0)}(K)$. The invertibility of $Q(L_\text{CME}-\Omega_*)Q$ follows for $\eps>0$ small enough from assumption (A2$_\text{CME}$). We do not employ this extended ansatz here as the leading order asymptotics suffice as an initial guess for our numerical computation of the solution $\vec{B}$ of \eqref{E:CME-stand}. 

Based on the above asymptotics one expects that under assumptions (A1$_\text{CME}$)-(A3$_\text{CME}$)  the ansatz $\vec{B}_\text{ans}$ with a localized solution $C:\R^d\to \C$ of \eqref{E:NLS-eff} provides an approximation of a localized solution of \eqref{E:CME-stand}. This can be made rigorous using a Lyapunov-Schmidt reduction argument. we refrain from this analysis here. In the following section we inspect whether assumptions (A1$_\text{CME}$)-(A3$_\text{CME}$) can be satisfied for $N=2,3,4$.

%---------------------------------------------------------------------
%---------------------------------------------------------------------
%---------------------------------------------------------------------
%CMEs for N=2,3,4, gaps, shape of the disp. rel.
\section{CMEs for $N$ Linearly Resonant Modes with $N=2,3,4$}
\label{S:CMEs_band-str}

We consider the perturbed PNLS \eqref{eq:PNLS-W} with $W$ in (A4). We study the dispersion relation $K\mapsto \Omega(K)$ of the resulting CMEs, where $K\in \R^d$ is the dual variable to $X$ with respect to the Fourier transform. A gap in the range of $\Omega$ is then a spectral gap of the corresponding spatial operator, i.e. of $L_{\text{CME}}(\nabla)$.
We concentrate on the existence of a spectral gap as one can typically expect the existence of standing solitary waves $\vec{A}(X,T)=e^{-\ri \Omega T}\vec{B}(X)$ with an (exponentially) localized $\vec{B}$, a solution of $\Omega \vec{B} - L_{\text{CME}}(\nabla)\vec{B} + \vec{N}(\vec{B})=0$, for $\Omega \notin \text{spec}(L_\text{CME}(\nabla))$. 

\subsection{Two modes: $N=2$}
Let us first choose $N=2$ in \eqref{eq:uapp}. The CMEs \eqref{E:CME} have the form
\beq
\begin{aligned}\label{eq:CME-2}
\ri\left(\pa_T + v_g^{(1)}\cdot \nabla\right)A_1 +\kappa_{11}A_1+ \kappa_{12} A_2 +(\gamma_1^{(1,1,1)}|A_1|^2+2\gamma_1^{(2,2,1)}|A_2|^2)A_1 + \delta \gamma_1^{(2,1,2)}A_2^2\overline{A}_1&=0\\
\ri\left(\pa_T + v_g^{(2)}\cdot \nabla\right)A_2 +\kappa_{22}A_2 + \overline{\kappa_{12}} A_1 +(\gamma_2^{(2,2,2)}|A_2|^2+2\gamma_2^{(1,1,2)}|A_1|^2)A_2 + \delta \gamma_2^{(1,2,1)}A_1^2\overline{A}_2&=0,
\end{aligned}
\eeq 
where
$$\delta:=\begin{cases} 0 \text{ if } k^{(1)}-k^{(2)}\notin \tfrac{1}{2}\Z^d \\ 1 \text{ if } k^{(1)}-k^{(2)}\in \tfrac{1}{2}\Z^d.\end{cases} $$
This is because only for $k^{(1)}-k^{(2)}\in \tfrac{1}{2}\Z^d$ is $2k^{(2)}-k^{(1)}\in k^{(1)}+\Z^d$ such that the term $A_2^2\overline{A_1}$ appears in the first equation. Similarly $A_1^2\overline{A_2}$ appears in the second equation only then. 
 
An example of a simple $W$ leading generically to $\kappa_{11},\kappa_{22},\kappa_{12}\neq 0$ is $W(x)=2\cos((k^{(1)}-k^{(2)})\cdot x)+1$. In that case is 
$$
\begin{aligned}
\kappa_{jj}&=-1-\delta \int_\T |p_j(x)|^2\cos(2(k^{(2)}-k^{(1)})\cdot x)\dd x, j=1,2\\
\kappa_{12}&=-\int_\T p_2(x)\overline{p_1}(x)\dd x - \delta \int_\T p_2(x)\overline{p_1}(x)e^{2\ri (k^{(2)}-k^{(1)})\cdot x }\dd x.
\end{aligned}
$$

In general, the dispersion relation $K \mapsto\Omega(K)$ of \eqref{eq:CME-2} is given by
\beq\label{eq:dr-2}
\text{det} \bspm \Omega-v_g^{(1)}\cdot K+\kappa_{11} & \kappa_{12} \\ \overline{\kappa_{12}} & \Omega -v_g^{(2)}\cdot K+\kappa_{22} \espm =0.
\eeq
In one dimension $d=1$ we have $v_g^{(2)}=-v_g^{(1)}$ and $p_2=\overline{p_1}$ due to the band structure symmetry $\omega(k)=\omega(-k)$ and the fact that at most two Bloch waves exist for each frequency $\omega_0$. Then $\kappa_{11}=\kappa_{22}\in \R$ and thus 
$$\Omega(K)=-\kappa_{11}\pm\sqrt{v_g^{(1)^2}K^2+|\kappa_{12}|^2},$$  
such that there is the spectral gap $(-\kappa_{11}-|\kappa_{12}|,-\kappa_{11}+|\kappa_{12}|)$. In the one dimensional case CMEs \eqref{eq:CME-2} with $v_g^{(2)}=-v_g^{(1)}$, $\kappa_{12}\in \R$, $\delta \gamma_1^{(2,1,2)}=\delta \gamma_2^{(1,2,1)}=0,$ and $\gamma_1^{(1,1,1)}=\gamma_1^{(2,2,1)}=\gamma_2^{(2,2,2)}=\gamma_2^{(1,1,2)}$ are well known to have a family of explicitly known solitary waves parametrized by the velocity $v\in (-v_g^{(1)},v_g^{(1)})$ and the frequency $\Omega \in (-\kappa_{11}-|\kappa_{12}|,-\kappa_{11}+|\kappa_{12}|)$, see \cite{AW89}. In the general case (with an spectral gap) solitary waves can be constructed numerically, see \cite{D14}.

For $d\geq 2$ we show that a gap exists only if $v_g^{(1)}$ and $v_g^{(2)}$ are linearly dependent and point in opposite directions. The non-existence of a gap is equivalent to the solvability of \eqref{eq:dr-2} in $K\in \R^d$ for all $\Omega \in \R$. Writing $K=\rho j$ with $\rho \in \R$ and $j\in \R^d, |j|=1$, we pose this as a problem in $\rho$ for fixed $\Omega\in \R$ and $j\in \R^d$. For the non-existence of a gap there must be for each $\Omega \in \R$ at least one direction $j$ such that a real solution $\rho$ exists. Assuming $j\cdot v_g^{(1)}\neq 0, j\cdot v_g^{(2)}\neq 0$, we get
\beq\label{eq:rho-2}
\begin{aligned}
\rho = \frac{1}{2(j\cdot v_g^{(1)})(j\cdot v_g^{(2)})}\left[\Omega j\cdot (v_g^{(1)}+v_g^{(2)})+\kappa_{11}j\cdot v_g^{(2)} +\kappa_{22}j\cdot v_g^{(1)}\right.\\
\left. \pm \left( (j\cdot v_g^{(1)}(\Omega+\kappa_{22})-j\cdot v_g^{(2)}(\Omega+\kappa_{11}) )^2 +4|\kappa_{12}|^2j\cdot v_g^{(1)} j\cdot v_g^{(2)}\right)^{1/2}\right].
\end{aligned}
\eeq 
We write $j\cdot v_g^{(l)}=|v_g^{(l)}|\cos(\theta_l)$, $l=1,2$ with $\theta_l\in [-\pi/2,\pi/2]$. A real $\rho$ exists if and only if 
\beq\label{eq:solv-cond-2}
\left(|v_g^{(1)}|\cos(\theta_1)(\Omega+\kappa_{22})-|v_g^{(2)}|\cos(\theta_2)(\Omega+\kappa_{11}) \right)^2 \geq -4|\kappa_{12}|^2 |v_g^{(1)}||v_g^{(2)}|\cos(\theta_1)\cos(\theta_2).
\eeq
There are two cases: $\cos(\theta_1)=-\cos(\theta_2)$ (i.e. $v_g^{(1)}=-\alpha v_g^{(2)}$ with $\alpha >0$) and $\cos(\theta_1)\neq -\cos(\theta_2)$. In the former case \eqref{eq:solv-cond-2} simplifies to 
\begin{align*}
 \left(\Omega+\frac{\kappa_{22}|v_g^{(1)}|+\kappa_{11}|v_g^{(2)}|}{|v_g^{(1)}|+|v_g^{(2)}|}\right)^2 \geq 4|\kappa_{12}|^2\frac{|v_g^{(1)}||v_g^{(2)}|}{(|v_g^{(1)}|+|v_g^{(2)}|)^2}.
\end{align*}
Hence, there is the spectral gap 
$$\left(-\frac{\kappa_{22}|v_g^{(1)}|+\kappa_{11}|v_g^{(2)}|}{|v_g^{(1)}|+|v_g^{(2)}|} - 2|\kappa_{12}|\frac{\sqrt{|v_g^{(1)}||v_g^{(2)}|}}{|v_g^{(1)}|+|v_g^{(2)}|},-\frac{\kappa_{22}|v_g^{(1)}|+\kappa_{11}|v_g^{(2)}|}{|v_g^{(1)}|+|v_g^{(2)}|} + 2|\kappa_{12}|\frac{\sqrt{|v_g^{(1)}||v_g^{(2)}|}}{|v_g^{(1)}|+|v_g^{(2)}|}\right).$$
Note that due to the continuous dependence of the solution $\rho$ on $j$, the excluded directions with $j\cdot v_g^{(1)}=0$ or $j\cdot v_g^{(2)} =0$ cannot change the existence of this gap.

\begin{figure}[ht!]
\includegraphics[scale=0.6]{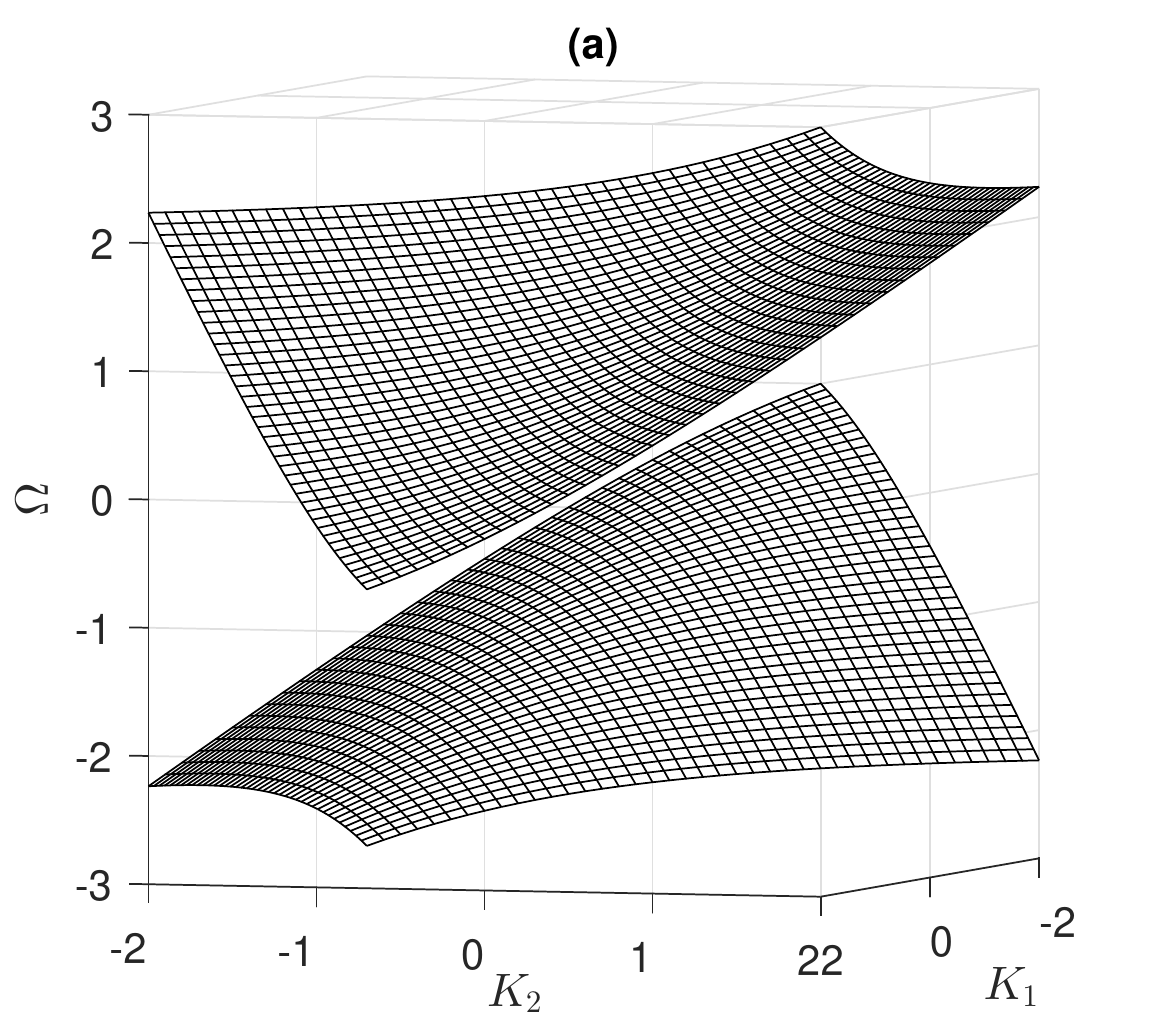}
\includegraphics[scale=0.6]{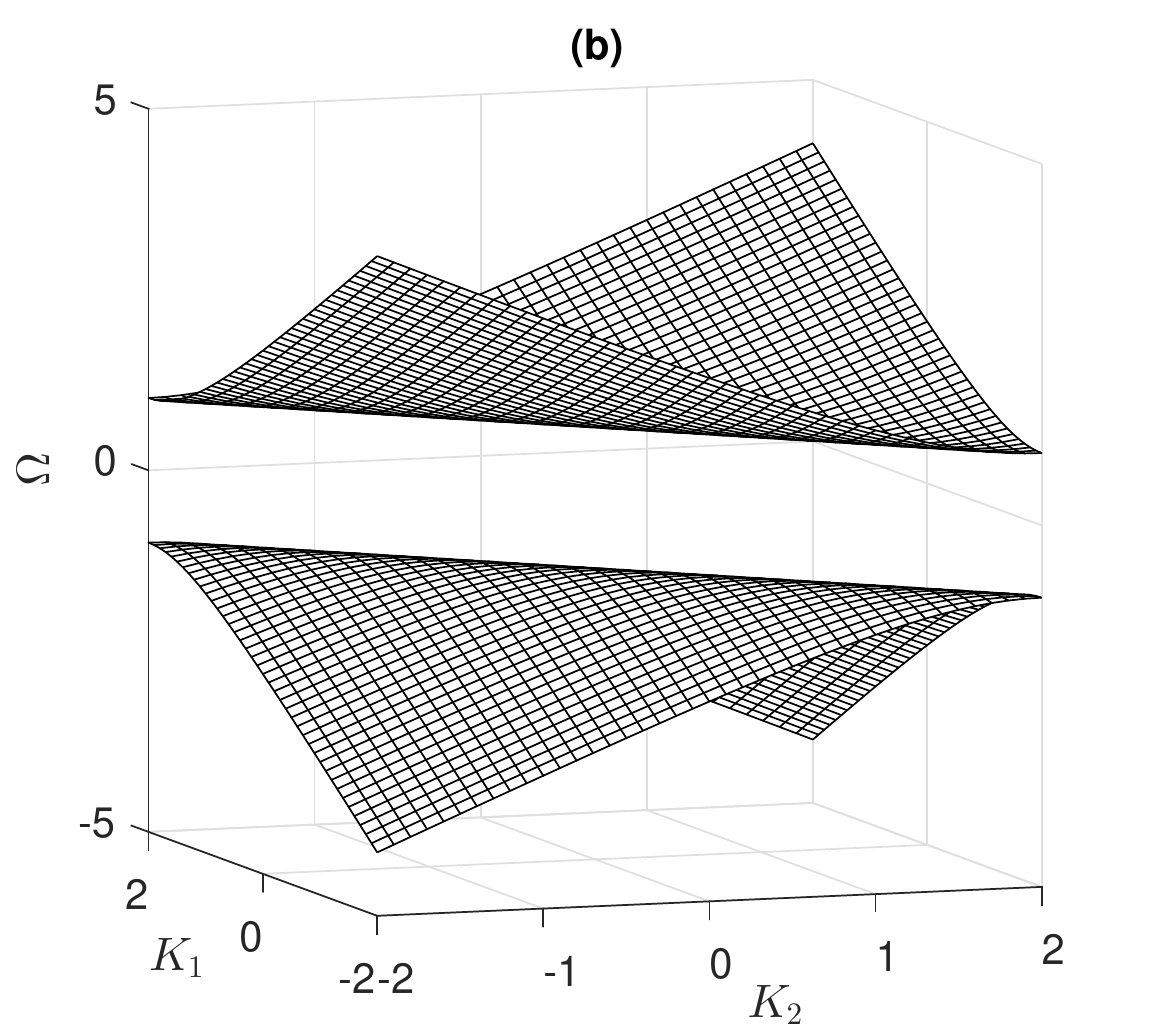}
\caption{Dispersion relation of \eqref{eq:CME-2} with $\kappa_{11}=\kappa_{22}=0, \kappa_{12}=\kappa_{21}=1$. $v_g^{(1)}=(0,1), v_g^{(2)} =(1,0)$ in (a) and $v_g^{(1)}=-v_g^{(2)} =(1,1)$ in (b).}
\label{F:N2-disprel}
\end{figure}

In the latter case ($\cos(\theta_1)\neq -\cos(\theta_2)$)  it is either $v_g^{(1)}=\alpha v_g^{(2)}$ with $\alpha >0$ or $v_g^{(1)}$ and $v_g^{(2)}$ are linearly independent. If $v_g^{(1)}=\alpha v_g^{(2)}, \alpha>0$, then $\cos(\theta_1)=\cos(\theta_2)$ and  the right hand side of \eqref{eq:solv-cond-2} is non-positive, such that \eqref{eq:solv-cond-2} always holds and there is  no spectral gap. If $v_g^{(1)}$ and $v_g^{(2)}$ are linearly independent (i.e. $|\cos (\theta_1)|\neq |\cos (\theta_2)|$), then one can always choose $j$ such that $\cos(\theta_1)$ and $\cos(\theta_2)$ have equal signs leading to a non-positive right hand side of \eqref{eq:solv-cond-2}. Hence there is at least one direction in which \eqref{eq:solv-cond-2} holds such that no spectral gap exists. Figure \ref{F:N2-disprel} (a) shows $\Omega=\Omega(K)$ for the case $d=2,\kappa_{11}=\kappa_{22}=0, \kappa_{12}=\kappa_{21}=1$, and $v_g^{(1)}=(0,1), v_g^{(2)} =(1,0)$.

The shape of the dispersion relation graph $\Omega=\Omega(K)$ near the edge of a spectral gap is important for the asymptotic construction of solutions in Sec. \ref{S:CME-asymp}. Namely, the construction is based on a Taylor expansion of the band structure near an isolated localized extremum.  For $N=2$ and  $v_g^{(1)}=-\alpha v_g^{(2)},\alpha>0$ there is no isolated extremum of $\Omega(K)$ at the spectral edges. This follows from \eqref{eq:rho-2}, where the square root vanishes at both edges and the values of $\rho(\theta_1)$ for $\theta_1\in [0,\pi/2)$ span a whole interval. From \eqref{eq:rho-2} one namely gets $\rho =\tfrac{1}{2\alpha |v_g^{(2)}|\cos(\theta_1)}\left[\tfrac{\alpha-1}{\alpha+1}(\kappa_{22}+\alpha \kappa_{11}\mp 2\sqrt{\alpha}|\kappa_{12}|)+\kappa_{11}-\alpha \kappa_{22}\right]$. This effectively one dimensional shape of the dispersion relation reflects the fact that with linearly dependent $v_g^{(1)}$ and $v_g^{(2)}$ the CMEs \eqref{E:CME-stand} both have the derivative in the same direction. Hence, the existence of solitary waves localized in all directions cannot be expected.

Figure \ref{F:N2-disprel} (b) shows $\Omega=\Omega(K)$ for the case $d=2, N=2,\kappa_{11}=\kappa_{22}=0, \kappa_{12}=\kappa_{21}=1$, and $v_g^{(1)}=-v_g^{(2)} =(1,1)$.

%--------------------------
\subsection{Odd number of modes: $N\in 2\N +1$}

First, we note that three (or more) linearly independent modes at one frequency $\omega_0$ are possible only for $d\geq 2$. 

We interpret the dispersion relation $\det(L_\text{CME}(\ri K)-\Omega I_{N\times N})=0$ as an algebraic equation in $\rho:=|K|$ for a given direction $K/|K|\in \R^d$ and a given $\Omega \in \R$. Because $L_{\text{CME}}(\ri K)$ is a Hermitian matrix for each $K$, the determinant of $L_\text{CME}(\ri K)-\Omega I_{N\times N}$ is real for each $K$. Hence, for odd $N$
the dispersion relation is an algebraic equation of odd degree with real coefficients such that it has a real root for each $\Omega \in \R$. This implies the absence of spectral gaps for $N\in 2\N +1$.

%--------------------------
\subsection{Four modes: $N=4$}\label{S:gap-N4}

We do not attempt to prove the existence of a gap for the general system \eqref{E:CME} with $N=4$. Instead we show that in two dimensions ($d=2$) under a simple condition on $\kappa$ a gap exists in the following symmetric case
\beq\label{E:N4-as}
\begin{aligned}
&v_g^{(1)}=-v_g^{(2)}=:v,v_g^{(3)}=-v_g^{(4)}=:w, \\
&\kappa_{12}=\kappa_{34}=:\alpha_1, \\
&\kappa_{14}=\kappa_{32}=:\alpha_2,\\
&\kappa_{13}=\kappa_{42}=:\alpha_3,\\
&\kappa_{jj}=0, j=1,\dots,4, 
\end{aligned}
\eeq
with $\alpha_j\in\C, j=1,2,3$. We also assume the linear independence of $v$ and $w$. 

A simple example for which these symmetries hold is
$$k^{(1)}=(l_1,l_2)\in  \B, k^{(2)}=-k^{(1)}, k^{(3)}=(-l_1,l_2), k^{(4)}=-k^{(3)}$$
with $l_1\neq 0, l_2\neq 0,$ $n_1=\dots = n_4=:n_*$, with a simple eigenvalue $\omega_{n_*}(k^{(j)})$ for each $j=1,\dots, 4$ , and with
$$W(x)=a\cos((k^{(2)}-k^{(1)})\cdot x)+b\cos((k^{(4)}-k^{(1)})\cdot x)+c\cos((k^{(3)}-k^{(1)})\cdot x),$$
$a,b,c  \in \R$.

\begin{figure}[ht!]
\includegraphics[scale=0.6]{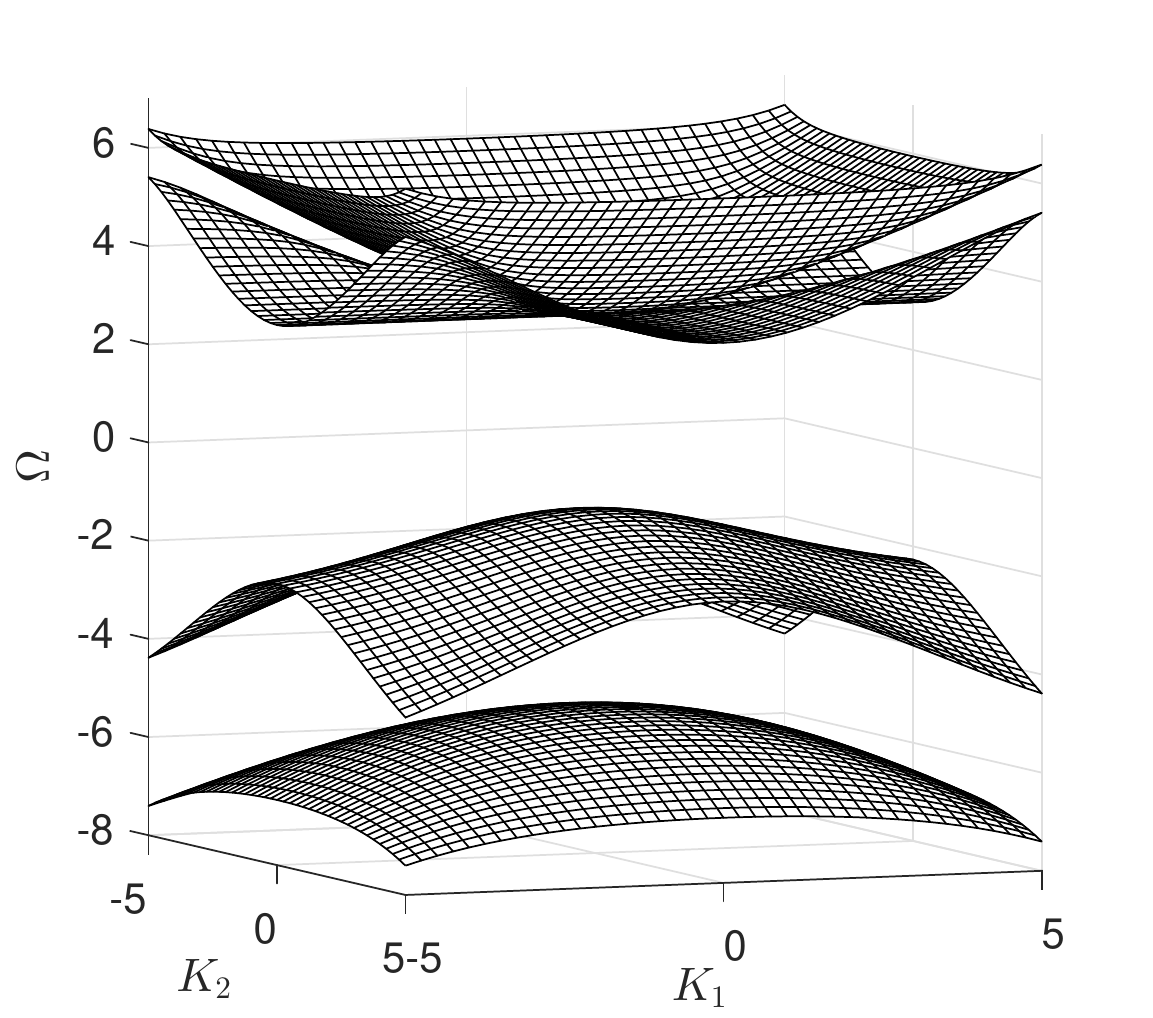}
\caption{Dispersion relation of \eqref{E:CME} with $N=4$, $v_g^{(j)}$ and $\kappa$ given in \eqref{E:N4-as}, and $\alpha_1 =3,\alpha_2=\alpha_3=1$, $v=(0,1), w=(1,0)$.}
\label{F:N4-disprel}
\end{figure}

After the transformation 
$$(\xi, \eta)^T=\phi(X_1,X_2):=\frac{1}{v_1 w_2-w_1v_2}\begin{pmatrix}w_2 X_1 - w_1 X_2\\ v_1 X_2 - v_2 X_1\end{pmatrix}$$
system \eqref{E:CME} becomes
\beq
\begin{aligned}\label{eq:CME-4}
\ri\left(\pa_T + \pa_\xi \right)A'_1 + \alpha_1 A'_2 + \alpha_3 A'_3 + \alpha_2 A'_4 +N_1(\vec{A'})&=0\\
\ri\left(\pa_T - \pa_\xi \right)A'_2 + \overline{\alpha_1} A'_1 + \overline{\alpha_2} A'_3 + \overline{\alpha_3} A'_4 +N_2(\vec{A'})&=0\\
\ri\left(\pa_T + \pa_\eta \right)A'_3 + \overline{\alpha_3} A'_1 + \alpha_2 A'_2 + \alpha_1 A'_4  +N_3(\vec{A'})&=0\\
\ri\left(\pa_T - \pa_\eta \right)A'_4 + \overline{\alpha_2} A'_1 + \alpha_3 A'_2 + \overline{\alpha_1} A'_3  +N_4(\vec{A'})&=0,\\
\end{aligned}
\eeq 
where $\vec{A}(\cdot,T)=\vec{A'}(\phi(\cdot),T)$. The dispersion relation of \eqref{eq:CME-4} is given by 
\beq\label{E:CME-DR}
\begin{aligned}
&\Omega^4-\Omega^2(K_\xi^2+K_\eta^2+2\sum_{j=1}^3|\alpha_j|^2)+4\Omega(\mbox{Re}(\overline{\alpha_1}\alpha_2\alpha_3)+\mbox{Re}(\alpha_1\overline{\alpha_2}\alpha_3))\\
&+(K_\xi K_\eta+|\alpha_2|^2-|\alpha_3|^2)^2+|\alpha_1|^2(K_\xi^2+K_\eta^2+|\alpha_1|^2)=2\mbox{Re}(|\alpha_1|^2\alpha_3^2+\alpha_2^2\overline{\alpha_1}^2).
\end{aligned}
\eeq
where $(K_\xi,K_\eta)$ is the dual variable to $(\xi,\eta)$. For $\Omega=0$ and $(K_\xi,K_\eta)\in \R^2$ the left hand side of this equation is larger than or equal to $|\alpha_1|^4$. The right hand side is smaller than or equal to $2|\alpha_1|^2(|\alpha_2|^2+|\alpha_3|^2)$. Hence, there is no real solution $(K_\xi,K_\eta)\in \R^2$ if $|\alpha_1|^2>2(|\alpha_2|^2+|\alpha_3|^2)$. As a result there is a spectral gap around 0 for \eqref{eq:CME-4} if $$|\alpha_1|^2>2(|\alpha_2|^2+|\alpha_3|^2).$$ 
For a slightly less general case of the matrix $\kappa$ this gap was found already in \cite{AP05}. Figure \ref{F:N4-disprel} shows an example of the dispersion relation of \eqref{E:CME} with $N=4$ and \eqref{E:N4-as} featuring a spectral gap.
The plot in Fig \ref{F:N4-disprel} suggests that assumptions (A1$_\text{CME}$)-(A3$_\text{CME}$) are satisfied at the lower edge of the spectral gap with $K_0:=(0,0)^T$. In Sec. \ref{S:num-stand} we compute corresponding solitary waves.

%-------------------------------------------------------------------
%-------------------------------------------------------------------
%-------------------------------------------------------------------

\section{Numerical Example of Standing Solitary Waves of the CMEs}
\label{S:num-stand}
In order to compute standing solitary waves $\vec{A}(X,t)=e^{-\ri \Omega T}\vec{B}(X)$, we perform a fixed point iteration on \eqref{E:CME-stand} starting with the initial guess given by \eqref{E:ansatzB} with a localized solution $C$ of \eqref{E:NLS-eff}. The fixed point iteration of \eqref{E:CME-stand} is performed using the Petviashvili iteration \cite{P76,AM05} in Fourier variables because standard finite difference or finite element discretizations of  \eqref{E:CME-stand} lead to spurious checkerboard oscillations caused by the decoupling of even and odd nodes- an effect referred to as Fermion doubling in the physics literature, see e.g. Sec. 5 in \cite{Hammer:2013}. 
%see also http://ilja-schmelzer.de/clm/doubling.php

To find a localized solution $C$ of \eqref{E:NLS-eff}, we first compute a radially symmetric $C(y)=C_\text{rad}(\rho)$ with $\rho:=|y|$ of \eqref{E:NLS-eff} with $D^2\eta^{(j_0)}$ replaced by $\mu I$, e.g. with $\mu:=\tfrac{1}{d}\sum_{j=1}^d \mu_j$, where $\{\mu_1,\dots,\mu_d\}$ are the eigenvalues of $D^2\Omega^{(j_0)}(K_0)$. We use the shooting method for the computation of $C_\text{rad}$. Then we perform a parameter continuation in order to deform $C_\text{rad}$ to a solution of \eqref{E:NLS-eff}. Hence, the algorithm is
\bi
\item[1)] shooting method for a solution $C_\text{rad}$ of 
\beq\label{E:Crad-eq}
\lambda C_\text{rad} + \mu(C_\text{rad}''(\rho) + \tfrac{d-1}{\rho}C_\text{rad}'(\rho)) +\Gamma C_\text{rad}^3(\rho)=0, \quad \rho \in (0,\infty),
\eeq
with $C_\text{rad}(\rho)\to 0$ as $\rho \to \infty$
\item[2)] parameter continuation in the parameter $\nu\in [0,1]$ (starting at $\nu=0$ and terminating at $\nu=1$) on equation \eqref{E:NLS-eff} with $D^2\eta^{(j_0)}$ replaced by $\mu I + \nu(D^2\eta^{(j_0)}-\mu I)$; finite difference discretization in the $x$-variable,
\item[3)] evaluation of \eqref{E:ansatzB} with $C$ from the previous step at $\nu=1$ and with a chosen (small) $\eps >0$,
\item[4)] Petviashvili iteration on \eqref{E:CME-stand} with $\Omega =\Omega_*+\eps^2\alpha$ using the initial guess given by \eqref{E:ansatzB},
\item[5)] parameter continuation in the parameter $\Omega$ from $\Omega_*+\eps^2\alpha$ to the desired value of $\Omega$ in the spectral gap of $L_\text{CME}$ (Petviashvili iteration at each continuation step).
\ei

\begin{example}\label{Ex:N4-gap}
We consider the two dimensional ($d=2$) case with $N=4$ and with the symmetries \eqref{E:N4-as}. The parameters are chosen like in Figure
\ref{F:N4-disprel}, i.e. $\alpha_1 =3,\alpha_2=\alpha_3=1$, $v=(0,1), w=(1,0)$. The spectral gap is approximately $(-1, 3)$. We let $\Omega_*$ be the lower edge, i.e. $\Omega_*\approx -1$. Using the notation of Sec. \ref{S:CME-asymp}, it is $\Omega_*=\Omega_{j_0}(K_0)$ with $j_0=2, K_0=0$, $\eta^{(j_0)}(K_0)=\tfrac{1}{2}(-1,-1,1,1)^T$. A numerical approximation via finite differences yields
$$D^2\Omega_{j_0}(K_0) \approx -\tfrac{1}{4}I_{2\times 2}, \Gamma \approx 2.125.$$
In Figure \ref{F:CME4-solit} (a) we plot the radial profile of the radially symmetric solution $C_\text{rad}$ of the NLS. Because $D^2\Omega_{j_0}(K_0)$ is a multiple of the identity, it is $C=C_\text{rad}$. Figure \ref{F:CME4-solit} (b) shows the first component of the asymptotic ansatz \eqref{E:ansatzB} with $\eps=0.1$. Due to $K_0=0$ and the realness of $C$ it is $\vec{B}_\text{ans}:\R^2\to \R^4$. Moreover, due to the form of $\eta^{(j_0)}(K_0)$ we have $B_{\text{ans},1}=B_{\text{ans},2}=-B_{\text{ans},3}=-B_{\text{ans},4}$. Figure \ref{F:CME4-solit} (c) shows the numerically computed (via the Petviashvili iteration) solution $\vec{B}$ at $\Omega = \Omega_*+\eps^2 = -0.99$. Only $B_1$ and $B_3$ are plotted since the computed solution satisfies $B_2=\overline{B_1}$, $B_4=\overline{B_3}$.  In Figure \ref{F:CME4-solit2} we plot the components $B_1$ and $B_3$ of the solution to \eqref{E:CME-stand} at $\Omega=-0.8$ obtained by a parameter continuation in $\Omega$ from $\Omega_*+\eps^2$. During the continuation the symmetry $B_1=B_2=-B_3=-B_4$ seems to break. Also note that the solution at $\Omega=-0.8$ is far from radially symmetric.

\begin{figure}[ht!]
\includegraphics[scale=0.6]{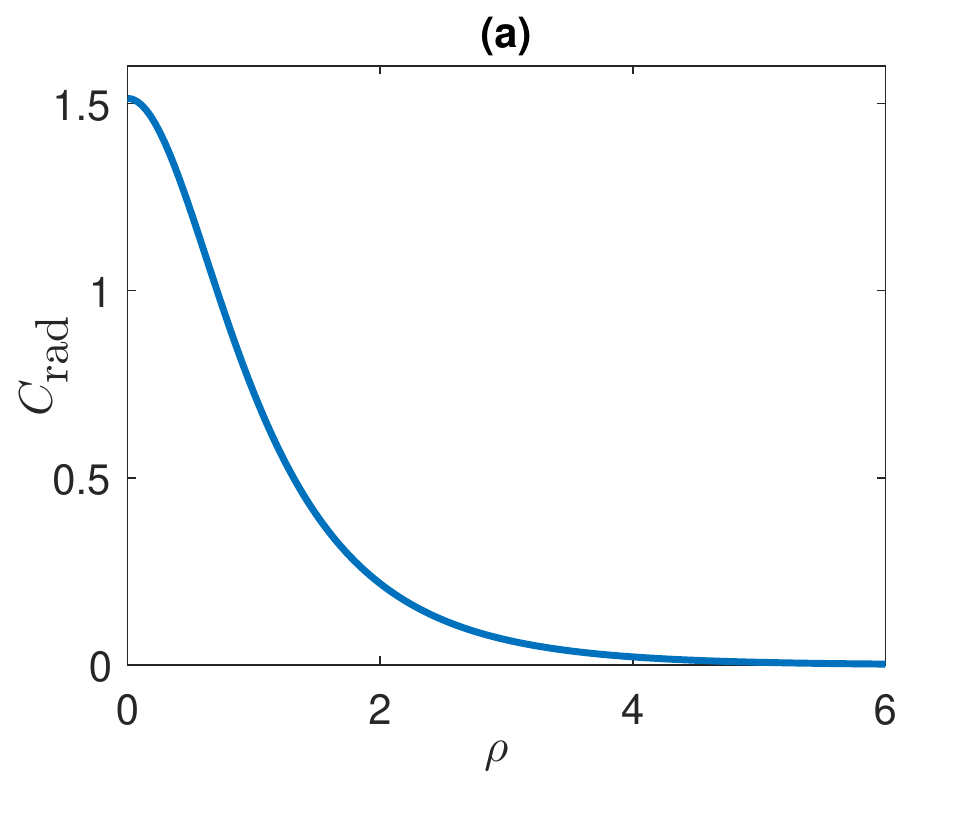}
\includegraphics[scale=0.6]{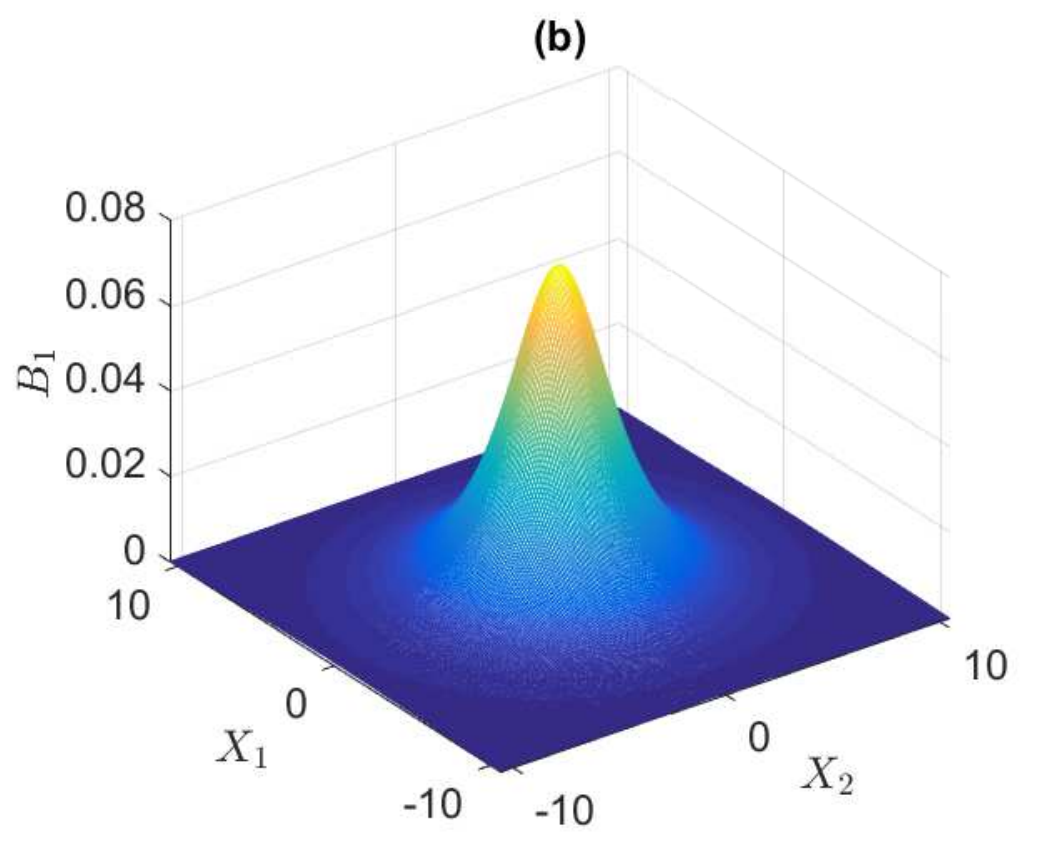}
\includegraphics[scale=0.6]{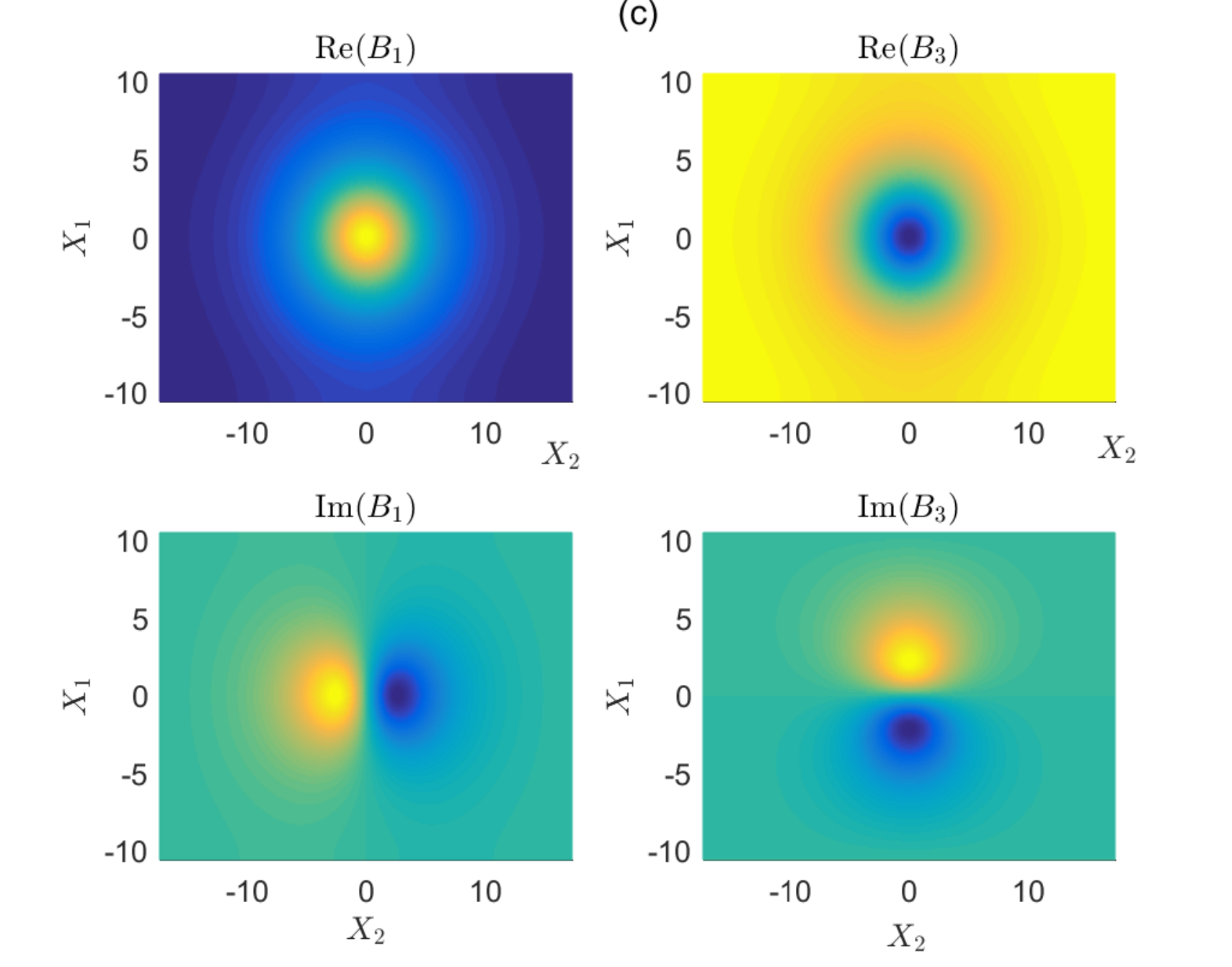}
\caption{Plots for Example \ref{Ex:N4-gap}. (a) The profile of the positive solution $C_\text{rad}$ of \eqref{E:Crad-eq}; (b) $B_{\text{ans},1}$ at $\eps =0.1$; (c) $B_1$ and $B_3$ at $\Omega = \Omega_*+\eps^2=-0.99$.}
\label{F:CME4-solit}
\end{figure}

\begin{figure}[ht!]
\includegraphics[scale=0.6]{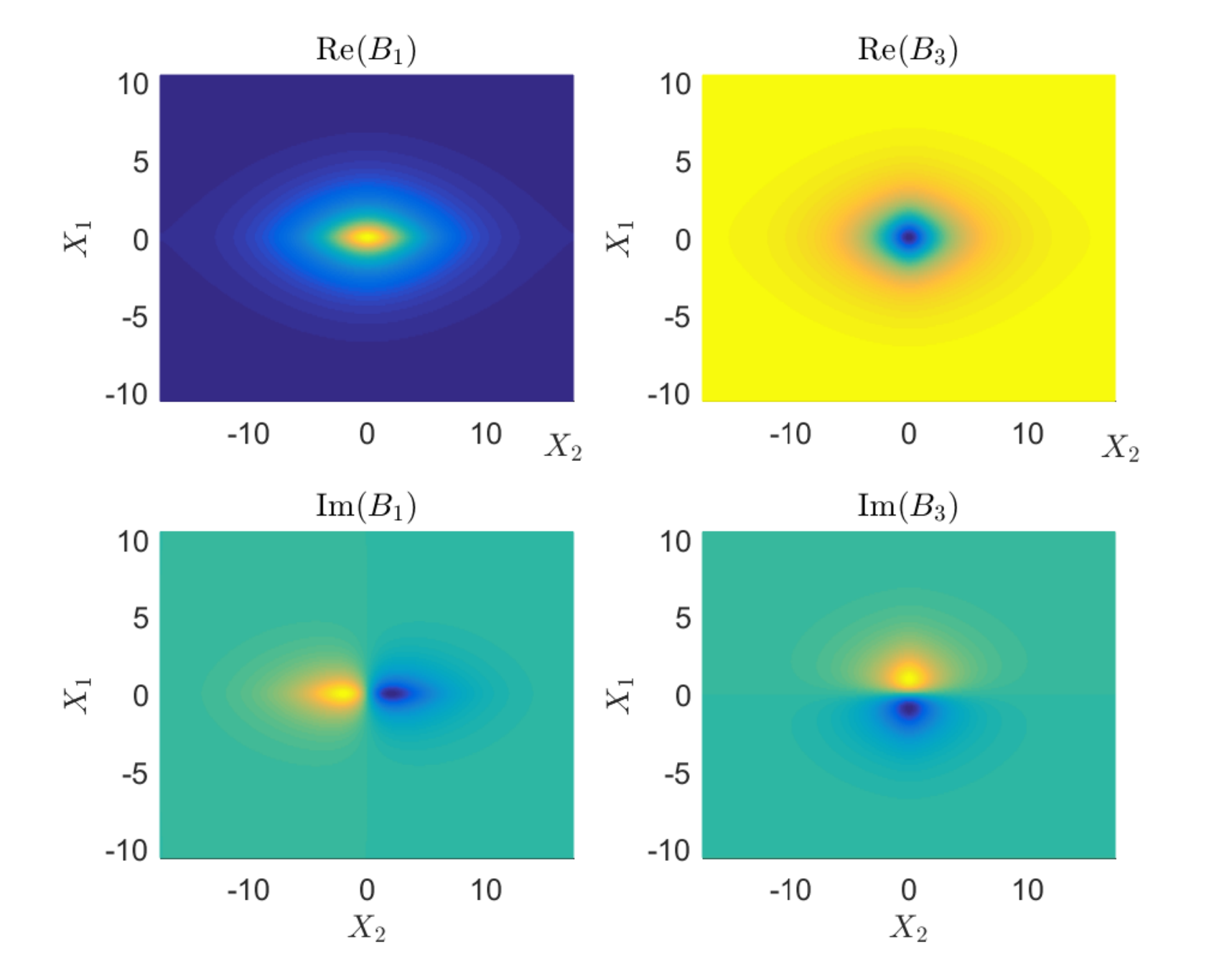}
\caption{Plots for Example \ref{Ex:N4-gap}. $B_1$ and $B_3$ at $\Omega = -0.8$.}
\label{F:CME4-solit2}
\end{figure}

\end{example}

%-------------------------------------------------------------------
%-------------------------------------------------------------------
%-------------------------------------------------------------------

\section{Moving Gap Solitons}\label{S:mov_solit}
An important feature of the CMEs in one dimension ($d=1$) is the existence of a whole family of solitary waves parametrized by their velocity. It is natural to ask whether this holds also for $d>1$. As explained in the introduction, the equal $\eps-$scaling of $t$ and $x$ in the asymptotic ansatz for a waves-packet of the PNLS would (in the affirmative case) imply that nearly solitary waves of the PNLS exist with a whole range of velocities. However, as we show next, for $N\leq 4$ and the most natural ansatz there seem to be no moving solitary waves in spectral gaps.

The simplest ansatz for a moving solitary wave of the time dependent CMEs \eqref{E:CME} is $\vec{A}(X,T)=e^{-\ri \Omega^{(v)} T} \vec{B}^{(v)}(X-vT)$ with some $v\in \R^d, \Omega^{(v)} \in \R$, and some localized $\vec{B}^{(v)}:\R^d \to \C^N$ ($|\vec{B}^{(v)}(Y)|\to 0$ as $|Y|\to\infty$). $\vec{B}^{(v)}$ must satisfy 
\beq\label{E:CME-mov}
\Omega^{(v)} \vec{B}^{(v)}=(L_\text{CME}(\nabla)+\ri v\cdot \nabla)\vec{B}^{(v)} - \vec{N}(\vec{B}^{(v)}).
\eeq
The corresponding dispersion relation is given by $K\mapsto \Omega^{(v)}(K),$ where 
$$\Omega_j^{(v)}(K)=\Omega_j^{(0)}(K)-v\cdot K, K\in\R^d, j\in \{1,\dots,N\}.$$
An exponentially localized $\vec{B}^{(v)}$ is expected in spectral gaps, i.e. for $\Omega^{(v)}\in \R\setminus \cup_{j=1}^N\Omega_j^{(v)}(\R^d)$. We have, however, the following
\blem
Let $v\in \R^d$. If for some $j_*\in \{1,\dots,N\}$ and some direction $\xi\in \R^d$,  such that $\xi^Tv\neq 0$, it is $\Omega_{j_*}^{(0)}(K)\to \mbox{const.}$ for $K=r\xi, \R \ni r\to \pm\infty$, then $\R=\Omega_{j_*}^{(v)}(\R\xi)$. Hence, in this case, there is no spectral gap of \eqref{E:CME-mov}.
\elem
\bpf
The graph $K\mapsto v\cdot K$ is a hyperplane in $\R^{d+1}$. Therefore, if the eigenvalue
$\Omega_{j_*}^{(0)},j_*\in\{1,\dots,N\}$ is asymptotically horizontal along the direction $\xi \in \R^d$ not orthogonal to $v$,
then $\Omega_{j_*}^{(v)}(K)=\Omega_j^{(0)}(K)-v\cdot K \to \infty$ for $K=r\xi, \R \ni r\to \infty$ or $ r\to -\infty$. In the opposite direction $\Omega_{j_*}^{(v)}(K) \to -\infty$. Due to  the intermediate value theorem it is $\Omega_{j_*}^{(v)}(\R\xi) = \R$.
\epf

A natural approach to constructing moving solitary waves is a parameter continuation from standing solitary waves (i.e. those with $v=0$).  Based on the discussion in Sec. \ref{S:CMEs_band-str} we do not expect (exponentially) localized solitary waves for $N\leq 3$ and $d\geq 2$. For $N=4$ and $d=2$ Example \ref{Ex:N4-gap} provides a setting with a spectral gap and standing solitary waves. Unfortunately, for $N=4$ and the symmetric case \eqref{E:N4-as} (satisfied by Example \ref{Ex:N4-gap}) the above mentioned horizontal asymptotics of the dispersion relation graph are always satisfied:
\blem
Let $N=4, d=2$ and let \eqref{E:N4-as} be satisfied. Then for each $\gamma \in \R\setminus (-|\alpha_1|,|\alpha_1|)$ 
there exists a direction $\theta \in (-\pi,\pi]$ and $j\in\{1,\dots,4\}$ such that $\Omega_j((r\cos(\theta),r\sin(\theta)))\to \gamma$ as $r\to\infty$.
\elem
\bpf
Writing $K_\xi=r \cos(\theta)$ and $K_\eta=r\sin(\theta)$ with $r>0, \theta \in (-\pi,\pi]$, the lemma is proved if the dispersion relation \eqref{E:CME-DR}
has a solution $\theta \in (-\pi,\pi]$ for each $|\Omega|\geq|\alpha_1|$ and for all $r>0$ large enough. 
Dispersion relation \eqref{E:CME-DR} can be rewritten as 
$$\left(\tfrac{r^2}{2}\sin(2\theta)+|\alpha_2|^2-|\alpha_3|^2\right)^2=r^2(\Omega^2- |\alpha_1|^2)-\psi,$$
where 
$$
\begin{aligned}\psi:=&\Omega^4-2\Omega^2\sum_{j=1}^3|\alpha_j|^2 + 4\Omega(\mbox{Re}(\overline{\alpha_1}\alpha_2\alpha_3+\alpha_1\overline{\alpha_2}\alpha_3))+|\alpha_1|^4-2|\alpha_1|^2\mbox{Re}(\alpha_3^2)-2\mbox{Re}(\overline{\alpha_1}^2\alpha_2^2).
\end{aligned}
$$
Hence
$$\sin(2\theta)=\pm \frac{2}{r^2}\left(|\alpha_3|^2-|\alpha_2|^2 \pm \sqrt{r^2(\Omega^2- |\alpha_1|^2)-\psi}\right).$$
For $r$ large a real solution $\theta$ clearly exists if and only if $|\Omega|\geq |\alpha_1|.$
%see disp_rel_CME_horiz_asymp_transformed.nb
\epf

\begin{figure}[ht!]
\includegraphics[scale=0.6]{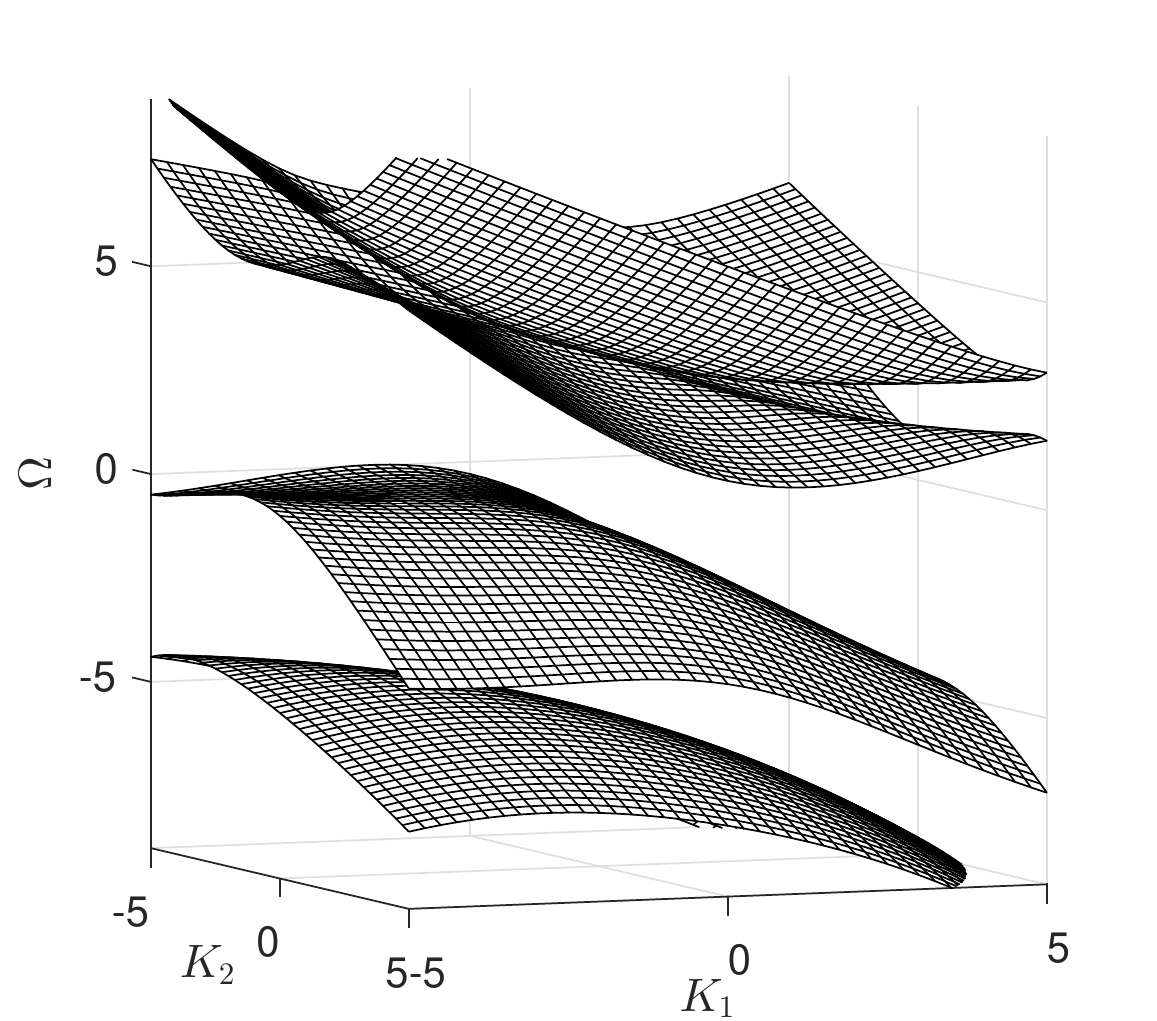}
\caption{Dispersion relation for \eqref{E:CME-mov} with $d=2,N=4$, $v=(0.3,0.3)^T$ and $v_g^{(j)},j=1,\dots,4,$ and $\kappa$ like in Example \ref{Ex:N4-gap}.}
\label{F:CME4mov-nogap}
\end{figure}

Figure \ref{F:N4-disprel} shows that for Example \ref{Ex:N4-gap} the eigenvalue $\Omega^{(0)}_3(K)$ converges to the constant $\alpha_1=3$ (the upper gap edge) along two lines in $\R^2$. In Figure \ref{F:CME4mov-nogap} we plot the dispersion relation $(K,\Omega^{(v)}(K))$ for the velocity $v=(1,1)^T$. Clearly, no gap occurs.

%-------------------------------------------------------------------
%-------------------------------------------------------------------

%Justification of CME in R^d
\section{Justification of CMEs \eqref{E:CME} as Amplitude Equations}
\label{S:justif}

In this section we prove Theorem \ref{T:justif}. This provides a rigorous justification of the CMEs \eqref{E:CME} as amplitude equations for waves-packets approximated by the formal ansatz \eqref{eq:uapp}.

In theorem \ref{T:justif} the functions $\hat{A}_j$ are Fourier transformations of $A_j$, where 
$$\hat{f}(k):=\frac{1}{(2\pi)^d}\int_{\R^d}f(x)e^{-\ri k\cdot x}\dd x$$
with the inverse $f(x)=\int_{\R^d}\hat{f}(k)e^{\ri k\cdot x}\dd k$. The space $L^1_r(\R^d)$ is
\begin{alignat}{2}
L^1_r\left(\R^d\right)\coloneqq \left\{ f \in L^1\left(\R^d\right) : \|f\|_{L^1_r(\R^d)} \coloneqq \int_{\R^d} \left( 1 + |x| \right)^r |f(x)| \dd x < \infty \right\}. 
\end{alignat}

We first summarize the functional analytic tools and some well known results regarding these. A more detailed presentation can be found, e.g., in \cite{DR18}.

We use the Bloch transformation 
\beq 
\begin{aligned}
&\cT: \quad f\mapsto \cT (f):=\widetilde{f},\\
&\widetilde{f}(x,k):=\sum_{\eta \in\Z^d}\widehat{f}\left(k+\eta \right)e^{\ri \eta \cdot x},
\end{aligned}
\eeq
which is an isomorphism from $H^s(\R^d)$ to $L^2(\B,H^s(\T))$ for $s\geq 0$ with the inverse 
$$f(x)=\int_\B \widetilde{f}(x,k)e^{\ri k\cdot x}\dd k.$$ 
The following properties hold for all $f,g\in H^s(\R^d)$ and $1\leq j \leq d$, $\N \ni p \leq s$:
\beq\label{E:perinxandk}
\cT (f)\left(x,k+e_j\right)=e^{-\ri x_j}\cT(f)(x,k),
\eeq
\beq\label{E:xandtder}
\cT (\partial_{x_j}^pf)(x,k)=(\partial_{x_j}+\ri k_j)^p(\cT f)(x,k),
\eeq
\beq\label{E:commper}
\cT(Vf)(x,k)=V(x)\cT(f)(x,k), \text{ if } V\in C(\T,\C),
\eeq
\beq\label{E:conv}
\cT(fg)(x,k)=((\cT f)*_{\B} (\cT g))(x,k):=\int_{\B}(\cT f)(x,k-l)(\cT g)(x,l)\mathrm{d}l.
\eeq

It is, however, inconvenient to carry out the estimates of the residual and of the asymptotic error in the $H^s(\R^d)$-norm (or equivalently in the $L^2(\B,H^s(\T))$-norm in the Bloch variables). The reason is the loss of $\eps$-powers when evaluating the $L^2$-norm if a function of the form $f(\eps \cdot)$. We prove Theorem \ref{T:justif} by working in an $L^1$ space in the Bloch variables. Estimate \eqref{E:est-main} is then obtained by an analog of the Lebesgue-Riemann lemma. Namely, if $\tilde{f} \in L^1(\B,H^s(\T))$ with $s>d/2$ and such that \eqref{E:perinxandk} holds, then the function $f(x):=\int_\B \tilde{f}(x,k)e^{\ri k\cdot x}\dd k$ satisfies 
$$\sup_{x\in \R^d}|f(x)|\leq c\|\tilde{f}\|_{L^1(\B,H^s(\T))} \quad \text{and} \ u(x)\to 0 \text{ as } |x| \to \infty.$$
See Lemma 1 in \cite{DR18}. This $L^1$ approach is one of the major differences from \cite{GMS08}, where the authors work in the space $H^s_\eps$, see Sec. \ref{S:intro}. Note that our assumption (A3) corresponds to the assumption of a ``closed mode system of order 3'' in \cite{GMS08}. There, however, the error is estimated by $C\eps$ in the $H^s_\eps$-norm, which produces the estimate $C\eps^{1-d/2}$ in the $L^\infty(\R^d)$-norm. For $C\eps^{3/2}$, like in Theorem \ref{T:justif}, the assumption of a closed mode system of a higher order (i.e. a stricter assumption than (A3)) is needed in \cite{GMS08}. As mentioned above, our approach also prevents us from having to consider high order correction terms in the asymptotic expansion of the solution.

Next, because the eigenfunctions  $(p_n(\cdot,k))_{n\in \N}$ are complete in $L^2(\T)$, we define the expansion/diagonalization operator \bigskip
\beq\label{E:Diag}
\mathcal{D}:L^1\left(\B,H^s(\T)\right)\rightarrow \cX(s):=L^1(\B,l^2_{s/d}),\quad \widetilde{u}\mapsto \vec{U}:=(U_n)_{n\in\N},\eeq
where $\util(x,k)=\sum_{n \in \N}U_n(k) p_n(x,k), U_n(k)=\langle \widetilde{u}(\cdot,k,t),p_n(\cdot,k)\rangle_{L^2(\T)}.$
$\cD$ is an isomorphism, see Lemma 6 in \cite{DR18}. Note that
$\|\vec{U}\|_{\cX(s)}:=\int_{\B}\|\vec{U}(k)\|_{l^2_{s/d}}\mathrm{d}k$
and
\beq \label{E:l2sd}
l^2_{s/d}:=\left\lbrace\vec{v}=(v_n)_{n\in\N}\in l^2(\mathbb{R}^d):\|\vec{v}\|^2_{l^2_{s/d}}=\sum_{n\in\N}n^{\frac{2s}{d}}|v_n|^2<\infty\right\rbrace.
\eeq

%-------------------------------------------

\subsection{Choice of an extended ansatz}
We return now to the formal asymptotic ansatz \eqref{eq:uapp}. Applying $\cT$, we get
\begin{alignat}{2}
\tilde{u}^{\text{app}}(x,k,t)=\eps^{1/2-d} \sum\limits_{\eta \in \Z^d} \sum \limits_{j=1}^N 
 \hat{A}_j\left( \frac{k-k^{(j)}+\eta}{\eps}, \eps t \right)p_j(x) e^{-\ri \omega _0 t}e^{ \ri \eta \cdot x}. \label{eq:uapp_tilde}
\end{alignat}
The Gross-Pitaevskii equation \eqref{eq:PNLS-W} transforms to
\begin{alignat}{2}
\left( \ri \partial_t -\mathcal{L}(x,k)\right)\tilde{u}(x,k,t) -\eps \sum\limits_{m=-m_*}^{m_*} a_m \tilde{u}(x, k - l^{(m)},t) - \sigma(x) (\tilde{u}*_\B\tilde{\bar{u}}*_\B\tilde{u})(x,k,t)=0. 
 \label{eq:PNLS-transf}
\end{alignat}
Here we have used the fact that for $g(x):=f(x)e^{\ri l\cdot x}$ is $(\cT g)(x,k)=(\cT f)(x,k-l)$. \\[1em]

The Bloch transformed residual of $u^\text{app}$ in the Gross-Pitaevskii (GP) equation is
\beq\label{E:uapp-res}
\begin{aligned}
&\cT\big(\text{GP}(u^{\text{app}})\big)(x,k,t)= \eps^{3/2-d} e^{-\ri \omega_0 t} \sum \limits_{\eta \in \Z^d} e^{\ri \eta \cdot x} \sum \limits_{j=1}^N \left[\ri \partial_T \hat{A}_j\left(\frac{k-k^{(j)}+\eta}{\eps},T\right) p_j(x) \right.\\
&+ 2\ri \hat{A}_j\left(\frac{k-k^{(j)}+\eta}{\eps},T\right) \frac{k-k^{(j)}+\eta}{\eps} \cdot \nabla p_j(x)  - \sum \limits^{m_*}_{m=-m_*} a_m \hat{A}_j\left(\frac{k-k^{(j)}- l^{(m)}+\eta}{\eps},T\right) p_j(x) \\
& -\left. \eps \left|\frac{k-k^{(j)}+\eta}{\eps}\right|^2 \hat{A}_j\left(\frac{k-k^{(j)}+\eta}{\eps},T\right) p_j(x) -I(x,k,t) \right],
\end{aligned}
\eeq
where
\begin{align*}
I(x,k,t)&=\sigma(x)\sum \limits_{\substack{\alpha, \beta, \gamma \\ \eta^{(\alpha)}, \eta^{(\beta)}, \eta^{(\gamma)} } }p_\alpha  ~\bar{p}_\beta ~p_\gamma e^{\ri (\eta^{(\alpha)} -\eta^{(\beta)} +\eta^{(\gamma)}) \cdot x} \int \limits_{\frac{\B+k^{(\beta)}-k^{(\gamma)}-\eta^{(\beta)}+\eta^{(\gamma)}}{\eps}} 
\int \limits_{\frac{\B-k^{(\gamma)}+\eta^{(\gamma)}}{\eps}}\\
& \hat{A}_\alpha\left(\frac{k-k^{(\alpha)}+k^{(\beta)}-k^{(\gamma)} +\eta^{(\alpha)}-\eta^{(\beta)}+\eta^{(\gamma)}}{\eps}- \tilde{h},T \right) \hat{\bar{A}}_\beta\left(\tilde{h} -\tilde{l},T\right)
 \hat{A}_\gamma\left( \tilde{l},T\right) ~d\tilde{l} ~d\tilde{h}.
\end{align*}
The $O(\eps^{1/2-d}$) terms in $\cT\big(\text{GP}(u^{\text{app}})\big)$ vanish due to the fact that each $p_j$ is an eigenfunction of $\cL(\cdot,k^{(j)})$ with the same eigenvalue $\omega_0$. The nonlinearity term was calculated using straightforward variable transformations in the convolution integrals.
  
Note that for $\hat{A}_j(K,T),j=1,\dots, N$ concentrated near $K=0$ one can formally approximate $\hat{A}_j(K,T)$ by $\chi_{B_{\eps^{-1/2}}(0)}(K)\hat{A}_j(K,T)$. We show below in Sec. \ref{S:resest} that $I(x,\cdot,t)$ can then be approximated near $k=k^{(j)}$ (precisely on $B_{\eps^{1/2}}(k^{(j)})$) by
$$
\begin{aligned}
\sum \limits_{\substack{\alpha, \beta, \gamma\in \{1,\dots,N\}\\ k^{(\alpha)}-k^{(\beta)}+k^{(\gamma)}\in k^{(j)}+\Z^d}}  (\hat{A}_\alpha*\hat{\bar{A}}_\beta*\hat{A}_\gamma)\left(\frac{k-k^{(j)}}{\eps},T \right) e^{\ri(k^{(\alpha)}-k^{(\beta)}+k^{(\gamma)}-k^{(j)})\cdot x}
p_\alpha(x)  ~\bar{p}_\beta(x) ~p_\gamma(x).
\end{aligned}
$$ 
The residual \eqref{E:uapp-res} is concentrated (in $k$) near the points $k^{(j)}, j=1,\dots,N$ but also near other points in $\B$. These are, firstly, the concentration points of those terms in $I$, for which $k^{(\alpha)}-k^{(\beta)}+k^{(\gamma)}\notin \{ k^{(1)},\dots, k^{(N)}\}+\Z^d$. Secondly, the multiplication of $\uapp$ by $W$ generates terms with new concentration points, namely when $l^{(m)} + k^{(j)} \notin \{ k^{(1)},\dots, k^{(N)}\}+\Z^d$ for some $m\in\{-m_*,\dots,m_*\}$, $j\in \{1,\dots,N\}$.  These terms in $I$ and in $Wu^\text{app}$ appear at the leading order and need to be eliminated in order to achieve a small residual. For this we define
\beq\label{E:J}
J:= J_0\cup J_{N} \cup J_{W},
\eeq 
where
\begin{align*}
J_0\coloneqq & \{k^{(1)},\dots,k^{(N)}\},\\
J_{N} \coloneqq &\left[ \left( \left\{ k^{(\alpha)} - k^{(\beta)} + k^{(\gamma)} | (\alpha,\beta,\gamma) \in \{1,\dots,N\}^3 \right\}+\Z^d\right) \cap \B \right] \setminus J_0, \\
J_{W}  \coloneqq &\left[ \left( \left\{ k^{(j)} +l^{(m)} |, j \in \{1,\dots,N\}, m\in\{-m_*,\dots,m_*\}\right\} +\Z^d\right) \cap \B \right] \setminus J_0.
\end{align*}
We also define
$$I_0:=\{(n_1,k^{(1)}),\dots,(n_N,k^{(N)})\}.$$
Clearly, $J$ is finite. The set $J_N\cap J_W$ consists of the new concentration points discussed above. Below we modify the formal ansatz $u^\text{app}$ to an extended ansatz $u^\text{ext}$ by including correction terms concentrated near points in $J_N\cup J_W$. This is done in the $\vec{U}$-variables, i.e. after the application of $\cD$.

The equation for $\vec{U}$ is
\begin{alignat}{2}
\big( \ri \partial_t - \Omega(k) \big) \Uvec(k,t) + \eps \sum\limits_{m=-m_*}^{m_*} M^{(m)}(k) \Uvec (k - l^{(m)},t) +\vec{F}(\Uvec,\Uvec,\Uvec)(k,t)=0, \ k \in \B, t>0,\label{eq:DT-PNLS}
\end{alignat}
where for all $i,j\in \N$
\begin{alignat}{2}
&\Omega_{jj}(k)=\omega_j(k), \ \Omega_{ij}=0 \text{ if } i \neq j,\\
&M^{(m)}_{ij}(k)= -a_m \langle p_j(\cdot,k- l^{(m)}),p_i(\cdot,k)\rangle_{\T}, \text{ and }\\
&F_j(\Uvec,\Uvec,\Uvec)(k):=-\langle \sigma(\cdot) (\tilde{u}*_\B\tilde{\bar{u}}*_\B\tilde{u})(\cdot,k,t),p_j(\cdot,k)\rangle_{\T},
\end{alignat}
and where $ \tilde{u} (x,k,t)= \sum\limits_{n \in \mathbb{N}} U_n(k,t) p_n(x,k).$

As advertised above, our extended ansatz must contain terms supported near $l\in J_0$ as well as near every $l\in J_N\cup J_W$. We choose these $k-$supports compact with the radius $\eps^{1/2}$ for the terms concentrated near $k^{(1)},\dots,k^{(N)}$ and consequently with the radius $3\eps^{1/2}$ for $J_N$.  For simplicity the terms concentrated near $l\in J_W$ are chosen also with the support of radius $3\eps^{1/2}$ because $J_W\cap J_N$ is not necessarily empty. Hence, we define the compactly supported envelopes $\tilde{A}_j, j \in \{1,\dots,N\}$ and $\tilde{A}_{n,l}, (n,l)\in (\N\times J) \setminus \{(n_1,k^{(1)}),\dots,(n_N,k^{(N)})\}$. The supports are
\begin{alignat}{2}
&\text{supp}\left( \tilde{A}_j (\cdot,T)\right) \subset B_{\eps^{-1/2}}, j=1, \dots ,N,\\
&\text{supp}\left( \tilde{A}_{n,l} (\cdot,T)\right) \subset B_{3\eps^{-1/2}}, (n,l)\in (\N \times J) \setminus \{(n_j,k^{(j)}):j=1,\dots,N\},
\end{alignat}
where $B_r:=B_r(0)\subset \R^d$ is the ball of radius $r>0$ centered at zero.

This leads us to the following extended ansatz. For $k\in \R^d, t \in \R$ we set
\begin{gather}
\begin{aligned}
U_{n_j}^{\text{ext}}(k,t)=& \eps^{1/2-d}\sum \limits_{\eta \in \Z^d}\left[\sum_{\stackrel{m\in \{1,\dots,N\}}{n_m=n_j}} \tilde{A}_m \left( \frac{k-k^{(m)}+\eta}{\eps} , T\right) + \eps\sum_{\stackrel{l\in J}{(n_j,l)\notin I_0}}\tilde{A}_{{n_j},l}\left(\frac{k-l+\eta}{\eps}, T \right) \right]e^{-\ri \omega_0 t}, \\
& \quad \text{for } j=1,\dots,N, \\
U_{n}^{\text{ext}}(k,t)=& \eps^{3/2-d}\sum \limits_{\eta \in \Z^d} \sum_{l\in J} \tilde{A}_{n,l}\left(\frac{k-l+\eta}{\eps}, T \right)e^{-\ri \omega_0 t}, \ \text{for } n \in \N\setminus I_0.
\end{aligned} \label{eq:uext}
\end{gather}
Note that the $\eta$-sums ensure the $\Z^d$-periodicity of $U_n^\text{ext}(\cdot,t)$ on $\R^d$. Due to the assumption on the supports of $\tilde{A}_j(\cdot,T)$ and $\tilde{A}_{n,l}(\cdot,T)$, the extended ansatz is supported in compact neighborhoods of the points in $J$ and their $\Z^d$-shifts.
 
For the estimate of the residual it is convenient to separate the leading order part of the ansatz. Hence, we set 
$$\vec{U}^{\text{ext},0}(k,t):=\eps^{1/2-d}e^{-\ri \omega_0 t}\sum \limits_{\eta \in \Z^d} \sum \limits_{j=1}^N \tilde{A}_j \left( \frac{k-k^{(j)}+\eta}{\eps} , T\right)e_{n_j} \text{ and } \vv{U}^{\text{ext},1}:=\Uvecext - \vv{U}^{\text{ext},0}.$$ 
Here $e_n$ is the $n-$th Euclidean unit vector in $\R^\N$. Note that $\vec{U}^{\text{ext},0}(\cdot,t)$ is supported in $\cup_{j=1}^N B_{\eps^{1/2}}(k^{(j)})+\Z^d$ and $\vec{U}^{\text{ext},1}(\cdot,t)$ is supported in $\cup_{l\in J} B_{3\eps^{1/2}}(l)+\Z^d$.

With the assumption on the supports of the envelopes the $\eta-$sums in \eqref{eq:uext} can be reduced to finite sums if we restrict $k$ to the periodicity cell $\B$, i.e. $k\in \B$. Let us namely define for any $k\in \B$ the set 
\begin{equation}
Z_k \coloneqq \left\{ \eta \in \Z^d : k - \eta \in \overline{\B}\right\}.
\end{equation} 
Then for each $\tilde{A}_j-$term, $j=1, \dots,N$,  the sum is over $\eta\in Z_{k^{(j)}}$ and for each $\tilde{A}_{n,l}-$term, $l\in J, n \in \N$, the sum is over $\eta\in Z_{l}$.

%-------------------------------------------
\subsection{Estimate of the residual}\label{S:resest}
The residual for the extended ansatz is
\begin{equation}
\vec{\text{Res}}(k,t):=\big( \ri \partial_t -\Omega(k) \big) \Uvecext(k,t) +\eps \sum\limits_{m=-m_*}^{m_*} M^{(m)}(k) \Uvecext (k - l^{(m)},t)   + \vec{F}(\Uvecext,\Uvecext,\Uvecext)(k,t).
\end{equation} 

The most complicated term is the nonlinearity $\vec{F}$. However, we need to understand its detailed structure only for the leading order part $\vec{F}(\vec{U}^{\text{ext},0},\vec{U}^{\text{ext},0},\vec{U}^{\text{ext},0})$ because this generates the nonlinearity in the CMEs. The rest of $\vec{F}$ will be simply estimated. We have
$$F_n(\vec{U}^{\text{ext},0},\vec{U}^{\text{ext},0},\vec{U}^{\text{ext},0})(k,t)=-\langle \sigma(\cdot) (\tilde{u}^{\text{ext},0}*_\B\tilde{\bar{u}}^{\text{ext},0}*_\B\tilde{u}^{\text{ext},0})(\cdot,k,t),p_n(\cdot,k)\rangle_{\T}, \ n \in \N,$$
where 
$$\tilde{u}^{\text{ext},0}(x,k,t):=\sum\limits_{n \in \mathbb{N}} U_n^{\text{ext},0}(k,t) p_n(x,k)=\eps^{1/2-d}\sum_{j=1}^N\sum_{\eta\in \Z^d}\tilde{A}_j\left(\frac{k-k^{(j)}+\eta}{\eps},T\right)p_{n_j}(x,k)e^{-\ri \omega_0 t}.$$ 
The Bloch transformation of a function $u(x)$ satisfies $\widetilde{\widebar{u}}(x,k)=\widebar{\widetilde{u}}(x,-k)$. We impose this condition on the ansatz $\tilde{u}^{\text{ext},0}$. This translates to the condition $\overline{\tilde{A}}(K,T)=\overline{\tilde{A}}(-K,T)$, which is satisfied, in particular, if $\tilde{A}$ is based on the Fourier transform of $A$, see \eqref{E_Atil-def}. We get
$$
\widetilde{\widebar{u}}^{\text{ext},0}(x,k,t)=\eps^{1/2-d}\sum_{j=1}^N\sum_{\eta\in \Z^d}\tilde{\bar{A}}_j\left(\frac{k+k^{(j)}-\eta}{\eps},T\right)\overline{p_{n_j}}(x,-k)e^{\ri \omega_0 t}.
$$
A direct calculation leads to
\begin{align*}
&F_n(\vec{U}^{\text{ext},0},\vec{U}^{\text{ext},0},\vec{U}^{\text{ext},0})(k,t) =\eps^{3/2-d}\sum_{l\in J_N\cup J_0} \sum_{(\alpha,\beta,\gamma)\in \Lambda_{l}}f_n^{\alpha,\beta,\gamma}(k,t)
\end{align*}
where 
$$
\Lambda_l:=\{(\alpha,\beta,\gamma)\in \{1,\dots,N\}^3:k^{(\alpha)}-k^{(\beta)}+k^{(\gamma)}\in l+\Z^d\}
$$
is the set of those indices $\alpha,\beta,\gamma$, for which the convolution term is concentrated at $l$, 
\begin{align*}
f_n^{\alpha,\beta,\gamma}(k,t)&:=e^{-\ri \omega_0 t}\int_{B_{2\eps^{-1/2}}}\int_{B_{\eps^{-1/2}}}\sum_{\eta\in \Z^d}\Atil_\alpha\left(\frac{k-(k^{(\alpha)}-k^{(\beta)}+k^{(\gamma)})+\eta-\eps r}{\eps},T\right)\Abartil_\beta(r-s,T)\times \\
&\times \Atil_\gamma(s,T) b^{(n)}_{\alpha,\beta,\gamma}(k-(k^{(\gamma)}-k^{(\beta)})-\eps r,k^{(\beta)}-\eps(r-s),k^{(\gamma)}+\eps s,k) \dd s\dd r \\
&\text{ for } k \in B_{3\eps^{1/2}}(k^{(\alpha)}-k^{(\beta)}+k^{(\gamma)}),\\
b^{(n)}_{\alpha,\beta,\gamma}(q,r,s,k)&:=-\langle \sigma(\cdot)p_{n_\alpha}(\cdot,q)\overline{p_{n_\beta}}(\cdot,r)p_{n_\gamma}(\cdot,s),p_n(\cdot,k)\rangle_\T, \text{ for } q,r,s,k\in \R^d.
\end{align*}

We consider the residual separately for $(n,k)$ with $n=n_j$ and $k$ near $k^{(j)}$ (or its integer shifts) with $j=1,\dots,N$ and for other $(n,k)\in \N\times \R^d$.

Firstly, for $k\in B_{3\eps^{1/2}}(k^{(j)}), j=1,\dots,N$ we have
\begin{align*}
\text{Res}&_{n_j}(k,t)= \eps^{3/2-d}\left[ \ri \partial_T  \tilde{A}_j\left( \frac{k-k^{(j)}}{\eps} ,T\right) +\eps^{-1} (\omega_{0} - \omega_{n_j}(k))  \tilde{A}_j \left( \frac{k-k^{(j)}}{\eps},T \right) \right. \\
& + \sum_{r=1}^N \sum_{\stackrel{m\in \{-m_*,\dots,m_*\}}{k^{(r)}+l^{(m)}\in k^{(j)}+\Z^d}} M^{(m)}_{n_j,n_r}(k) \tilde{A}_r \left( \frac{k-k^{(j)}}{\eps} ,T\right) \\
&-\sum_{(\alpha,\beta,\gamma)\in \Lambda_{k^{(j)}}}\int_{B_{2\eps^{-1/2}}}\int_{B_{\eps^{-1/2}}}\Atil_\alpha\left(\frac{k-k^{(j)}-\eps r}{\eps},T\right)\Abartil_\beta(r-s,T)\Atil_\gamma(s,T) \times \\
&\left. \times b^{(n_j)}_{\alpha,\beta,\gamma}(k-(k^{(\gamma)}-k^{(\beta)})-\eps r,k^{(\beta)}-\eps(r-s),k^{(\gamma)}+\eps s,k) \dd s\dd r\right] e^{-\ri\omega_0 t}+\hot.
\end{align*}
We can recover $\text{Res}_{n_j}(k,t)$ in the neighborhood of $k^{(j)}+\Z^d$ by the $\Z^d$ periodicity in $k$. 

Secondly, for $(n,k)\in (\N\times \R^d) \setminus \cup_{j=1}^N (n_j,B_{3\eps^{1/2}}(k^{(j)})+\Z^d)$
\begin{align*}
\text{Res}&_{n}(k,t)= \eps^{3/2-d}\left[ (\omega_{0} - \omega_{n}(k)) \sum_{l\in J, \eta\in \Z^d} \tilde{A}_{n,l} \left( \frac{k-l+\eta}{\eps},T \right) \right. \\
& + \sum_{j=1}^N\sum_{\stackrel{m=-m_*,\dots,m_*}{\eta\in \Z^d}} M^{(m)}_{n,n_j}(k) \tilde{A}_j \left( \frac{k-k^{(j)}-l^{(m)}+\eta}{\eps} ,T\right) \\
&\left.+\sum_{l\in J_N\cup J_0} \sum_{(\alpha,\beta,\gamma)\in \Lambda_{l}}f_n^{\alpha,\beta,\gamma}(k,t)\right] e^{-\ri\omega_0 t}+\hot.
\end{align*}
The higher order terms $\hot$ above come from the $T-$derivative of $ \vv{U}^{\text{ext,1}}$, from the application of the potential $W$ to $ \vv{U}^{\text{ext,1}}$ and from the nonlinearity $F_n$ with at least one of the three arguments being $ \vv{U}^{\text{ext,1}}$.

\underline{1) $\text{Res}_{n_j}(k,t)$ for $k\in B_{\eps^{1/2}}(k^{(j)}), j=1,\dots,N$}

In order to ensure the smallness of $\text{Res}_{n_j}$, we set 
\beq\label{E_Atil-def}
\Atil_j(K,T):=\chi_{B_{\eps^{-1/2}}}(K) \Ahat_j(K,T), j\in \{1,\dots,N\},
\eeq
where $(A_1,\dots,A_N)$ solves \eqref{E:CME}. This makes sense because, as we show next, the leading order part of $\text{Res}_{n_j}$ is approximated by the Fourier transform of the left hand side of the $j$-th equation in \eqref{E:CME}. Let us start with the linear terms. For $k\in B_{\eps^{1/2}}(k^{(j)})$ we have 
$$|\omega_{n_j}(k)-\omega_0-\eps\frac{k-k^{(j)}}{\eps}\cdot v_g^{(j)}|\leq c \eps^2 \left|\frac{k-k^{(j)}}{\eps}\right|^2 (1+h_\eps(k))$$
with $h_\eps(k) \to 0$ for $\eps \to 0$. Hence
\beq\label{E:omA-est}
\begin{aligned}
&\int_{\B+k^{(j)}}\left|\eps^{-1}(\omega_0-\omega_{n_j}(k^{(j)}))\Atil_j\left(\frac{k-k^{(j)}}{\eps},T\right)-\chi_j(k)\frac{k-k^{(j)}}{\eps}\cdot v_g^{(j)}\Ahat_j\left(\frac{k-k^{(j)}}{\eps},T\right)\right|\dd k \\
&\leq c\eps^{1+d}\|\Ahat_j(\cdot,T)\|_{L^1_2(B_{\eps^{-1/2}}(0))}\leq c\eps^{1+d}\|\Ahat_j(\cdot,T)\|_{L^1_2(\R^d)},
\end{aligned}
\eeq
where $$\chi_j:=\chi_{B_{\eps^{1/2}}(k^{(j)})}.$$

Next, we compare the linear coupling terms.  For $j,r\in \{1,\dots,N\}$ let $m\in  \{-m_*,\dots,m_*\}$ be such that $k^{(r)}+l^{(m)}\in k^{(j)}+\Z^d$. We rewrite in $M^{(m)}_{n_j,n_r}(k)$ the variables as
$$k-l^{(m)}=k^{(r)}+\eps \frac{k-k^{(j)}}{\eps} + k^{(j)}- k^{(r)}-l^{(m)}$$
and
$$k=k^{(j)}+\eps \frac{k-k^{(j)}}{\eps}.$$
We define 
$$K:=\tfrac{k-k^{(j)}}{\eps} \text{  for  } k\in B_{\eps^{1/2}}(k^{(j)}), \ j=1,\dots,N,$$
and write
$\kappa_{jr}=\sum_{\stackrel{m\in\{-m_*,\dots,m_*\}}{k^{(r)}+l^{(m)}\in k^{(j)}+\Z^d}} \kappa_{jr}^{(m)}$, where
$$\kappa_{jr}^{(m)}:=-a_{m}\int_\T e^{\ri(k^{(r)}+l^{(m)}-k^{(j)})\cdot x} p_r(x)\overline{p_j}(x)\dd x.$$

Then for $k\in B_{\eps^{1/2}}(k^{(j)})$ and $m\in \{-m_*,\dots, m_*\}$ such that $k^{(r)}+l^{(m)}\in k^{(j)}+\Z^d$ it is
\beq\label{E:M-kap_est}
\begin{aligned}
|M^{(m)}_{n_j,n_r}(k)-\kappa_{jr}^{(m)}|\leq &|a_m|\int_\T \left|e^{\ri(k^{(r)}-k^{(j)}+l^{(m)})\cdot x}p_{n_r}(x,k^{(r)}+\eps K)\overline{p_{n_j}}(x,k^{(j)}+\eps K)\right.\\
& \qquad \left. -e^{\ri(k^{(r)}+l^{(m)}-k^{(j)})\cdot x}p_{n_r}(x,k^{(r)})\overline{p_{n_j}}(x,k^{(j)})\right|\dd x\\
\leq &|a_m|\int_\T \left|p_{n_r}(x,k^{(r)}+\eps K)\right|\left|p_{n_j}(x,k^{(j)}+\eps K)-p_{n_j}(x,k^{(j)})\right|\dd x \\
+&|a_m|\int_\T \left|p_{n_j}(x,k^{(j)})\right|\left|p_{n_r}(x,k^{(r)}+\eps K)-p_{n_r}(x,k^{(r)})\right|\dd x\\
\leq & c\eps K
\end{aligned}
\eeq
due to (A6).

From \eqref{E:M-kap_est} we get
\beq\label{E:MA-est}
\left\|M^{(m)}_{n_j,n_r}(\cdot)\Atil_r\left(\frac{\cdot-k^{(j)}}{\eps},T\right)-\kappa_{j,r}^{(m)}\chi_j\Ahat_r\left(\frac{\cdot-k^{(j)}}{\eps},T\right)\right\|_{L^1(\B+k^{(j)})}\leq c\eps^{d+1}\|\Ahat_r(\cdot,T)\|_{L^1_1(\R^d)}.
\eeq
For the nonlinearity in the leading order part of $\text{Res}_{n_j}(k,t)$ with $k\in B_{\eps^{1/2}}(k^{(j)})$ we need to estimate
$$
\begin{aligned}
\Psi_j^{\alpha,\beta,\gamma}(k,T):=&\int_{B_{2\eps^{-1/2}}}\int_{B_{\eps^{-1/2}}}\Atil_\alpha\left(\frac{k-k^{(j)}-\eps r}{\eps},T\right)\Abartil_\beta(r-s,T)\Atil_\gamma(s,T) \times \\
&\times b^{(n_j)}_{\alpha,\beta,\gamma}(k-(k^{(\gamma)}-k^{(\beta)})-\eps r,k^{(\beta)}-\eps(r-s),k^{(\gamma)}+\eps s,k) \dd s\dd r\\
&-\chi_j\gamma_j^{(\alpha,\beta,\gamma)}(\Ahat_\alpha *\Ahatbar_\beta * \Ahat_\gamma)\left(\frac{k-k^{(j)}}{\eps},T\right)
\end{aligned}
$$
for each $(\alpha,\beta,\gamma)\in \Lambda_{k^{(j)}}$. We write
$$\Psi_j^{\alpha,\beta,\gamma}=I_1+I_2+I_3,$$
where (using $\gamma_j^{(\alpha,\beta,\gamma)}=b^{(n_j)}_{\alpha,\beta,\gamma}(k^{(\alpha)},k^{(\beta)},k^{(\gamma)},k^{(j)})$)
$$
\begin{aligned}
I_1:=&\int_{B_{2\eps^{-1/2}}}\int_{B_{\eps^{-1/2}}}\Atil_\alpha\left(\frac{k-k^{(j)}-\eps r}{\eps},T\right)\Abartil_\beta(r-s,T)\Atil_\gamma(s,T) \times \\
&\times \left[b^{(n_j)}_{\alpha,\beta,\gamma}(k-(k^{(\gamma)}-k^{(\beta)})-\eps r,k^{(\beta)}-\eps(r-s),k^{(\gamma)}+\eps s,k)\right.\\
&\quad\left. -b^{(n_j)}_{\alpha,\beta,\gamma}(k^{(\alpha)},k^{(\beta)},k^{(\gamma)},k^{(j)}) \right]\dd s\dd r,
\end{aligned}
$$
$$
\begin{aligned}
I_2:=\gamma_j^{(\alpha,\beta,\gamma)}&\left[(\Atil_\alpha *_{B_{2\eps^{-1/2}}} \Abartil_\beta *_{B_{\eps^{-1/2}}}\Atil_\gamma)\left(\frac{k-k^{(j)}-\eps r}{\eps},T\right)\right.\\
&\quad \left.-\chi_{B_{3\eps^{1/2}}(k^{(j)})}(k)(\Ahat_\alpha *\Ahatbar_\beta * \Ahat_\gamma)\left(\frac{k-k^{(j)}}{\eps},T\right)\right],
\end{aligned}
$$
and
$$
I_3:=\gamma_j^{(\alpha,\beta,\gamma)}\chi_{B_{3\eps^{1/2}}(k^{(j)})\setminus B_{\eps^{1/2}}(k^{(j)})}(k)(\Ahat_\alpha *\Ahatbar_\beta * \Ahat_\gamma)\left(\frac{k-k^{(j)}}{\eps},T\right).
$$
We estimate $\|I_j(\cdot,T)\|_{L^1(\B+k^{(j)})}, j=1,2,3$ analogously to \cite{DR18}. For $\|I_1(\cdot,T)\|_{L^1(\B+k^{(j)})}$ we use the multilinearity of $b^{(n_j)}$ and the Lipschitz continuity of the Bloch functions $p_n(x,\cdot)$ in (A6). $\|I_{2,3}(\cdot,T)\|_{L^1(\B+k^{(j)})}$ are small due to the decay of $\Ahat_j(\cdot,T)$, namely due to
$$\left\|(1-\chi_{B_{a\eps^{-1/2}}})\Ahat_j(\cdot,T)\|_{L^1(\R^d)}\leq c_a \eps^{\nu/2} \|\Ahat_j(\cdot,T)\right\|_{L^1_\nu(\R^d)}$$
for any $a,\nu >0$. One obtains
\beq\label{E:Ij-est}
\begin{aligned}
&\|I_1(\cdot,T)\|_{L^1(\B+k^{(j)})}\leq c \eps^{d+1}\|\Ahat_\alpha(\cdot,T)\|_{L^1_1(\R^d)}\|\Ahat_\beta(\cdot,T)\|_{L^1_1(\R^d)}\|\Ahat_\gamma(\cdot,T)\|_{L^1_1(\R^d)},\\
&\|I_{2,3}(\cdot,T)\|_{L^1(\B+k^{(j)})}\leq c\eps^{d+\nu/2} \|\Ahat_\alpha(\cdot,T)\|_{L^1_\nu(\R^d)}  \|\Ahat_\beta(\cdot,T)\|_{L^1_\nu(\R^d)} \|\Ahat_\gamma(\cdot,T)\|_{L^1_\nu(\R^d)}
\end{aligned}
\eeq
for any $\nu>0$.

Estimates \eqref{E:omA-est},\eqref{E:MA-est}, and \eqref{E:Ij-est} ensure that CMEs \eqref{E:CME} yield a small $\text{Res}_{n_j}(k,t)$, $j=1,\dots,N$ for $k$ in a vicinity of $k^{(j)}+\Z^d$. The rest of the leading order part of $\vec{\text{Res}}$ is made small by a suitable choice of the correction terms $\Atil_{n,l}, n\in \N,l\in J$.

\underline{2) Rest of $\vec{\text{Res}}$ at $O(\eps^{3/2})$}

For $(n,k)\in \N\times \R^d\setminus \cup_{j=1}^N(n_j,B_{3\eps^{1/2}}(k^{(j)})+\Z^d)$ the leading order ($O(\eps^{3/2})$) part of the residual  vanishes if we set
\beq\label{E:Atil_def}
\Atil_{n,l}\left(\frac{k-l}{\eps},T\right):=\frac{1}{\omega_0-\omega_n(k)}\left[\sum_{j=1}^N\sum_{\stackrel{m\in \{-m_*,\dots,m_*\}}{k^{(j)}+l^{(m)}\in l + \Z^d}}M^{(m)}_{n,n_j}\Atil_j\left(\frac{k-l}{\eps},T\right) - \sum_{(\alpha,\beta,\gamma)\in \Lambda_l}f_n^{\alpha,\beta,\gamma}(k,t) \right]
\eeq
for all $l\in J$ if $n\in \N\setminus \{n_1,\dots,n_N\}$, and all $l \in J$ such that $(n_j,l)\notin I_0$ if $n=n_j, j=1,\dots,N$. 

Note that in \eqref{E:Atil_def} is 
$$|\omega_0-\omega_n(k)|>c n^{2/d}, \quad n \in \N$$
with some $c>0$. For $n\in \{n_1,\dots, n_N\}$ this holds because $(n,l)\neq (n_j,k^{(j)})$ and $k\in B_{3\eps^{1/2}}(l)$, and the growth in $n$ is guaranteed by \eqref{E:as-distr}.

Combining now \eqref{E:omA-est},\eqref{E:MA-est}, and \eqref{E:Ij-est}, we have
\beq\label{E:Res-est1}
\begin{aligned}
\|\vec{\text{Res}}(\cdot,t)\|_{\cX(s)}\leq &\|\hot(\cdot,t)\|_{\cX(s)}+ c \eps^{\frac{5}{2}} \left(\sum_{\alpha\in \{1,\dots,N\}} \|\Ahat_\alpha(\cdot,T)\|_{L^1_2(\R^d)}\right.\\
& \left.+\sum_{\alpha,\beta,\gamma\in \{1,\dots,N\}} \|\Ahat_\alpha(\cdot,T)\|_{L^1_2(\R^d)}\|\Ahat_\beta(\cdot,T)\|_{L^1_2(\R^d)}\|\Ahat_\gamma(\cdot,T)\|_{L^1_2(\R^d)}
\right).
\end{aligned}
\eeq
Thus it remains to estimate the $\hot$. 

The $n$-th component $(\hot)_n(k,t)$ consists of a (finite) sum over $l$ of the terms 
$\eps^{5/2-d}\pa_T\Atil_{n,l}$$\left(\frac{k-l}{\eps},T\right)e^{-\ri \omega_0 t},$ and of $\eps \sum_{m\in \{-m_*,\dots, m_*\}} (M^{(m)}(k)\vv{U}^{\text{ext},1}(k-l^{(m)},t))_n$, $F_n(\vv{U}^{\text{ext},1},\vv{U}^{\text{ext},0},\vv{U}^{\text{ext},0})$, $F_n(\vv{U}^{\text{ext},0},\vv{U}^{\text{ext},1},\vv{U}^{\text{ext},0})$ as well as of nonlinear terms quadratic or cubic in $\vv{U}^{\text{ext},1}$.

For  $\eps^{5/2-d}\pa_T\Atil_{n,l}\left(\frac{k-l}{\eps},T\right)e^{-\ri \omega_0 t}$ we first note that
$$
\left|M_{n,n_j}^{(m)}\right|=\frac{|a_m|}{|\omega_n(k)|^q}\langle \cL^q(\cdot,k) p_{n_j}(\cdot, k-l^{(m)}), p_n(\cdot,k)\rangle_\T \leq cn^{-2q/d}\sup_{k\in \B}\|p_{n_j}(\cdot, k)\|_{H^{2q}(\T)}
$$
for any $q\in \N$. Similarly,
$$
\begin{aligned}
|b^{(n)}_{\alpha,\beta,\gamma}(o,r,s,k)|&=\frac{1}{|\omega_n(k)|^q}\langle \cL^q(\cdot,k) ( \sigma(\cdot)p_{n_\alpha}(\cdot,o)\overline{p_{n_\beta}}(\cdot,r)p_{n_\gamma}(\cdot,s)),p_n(\cdot,k)\rangle_\T\\
&\leq cn^{-2q/d}\sup_{k\in \B}\|p_{n_\alpha}(\cdot, k)\|_{H^{2q(\T)}}\sup_{k\in \B}\|p_{n_\beta}(\cdot, k)\|_{H^{2q(\T)}}\sup_{k\in \B}\|p_{n_\gamma}(\cdot, k)\|_{H^{2q(\T)}}
\end{aligned}
$$
for  all $o,r,s,k \in 2\B$ and $q\in \N\cap (\frac{d}{4}, \infty)$ if $\sigma\in H^{2q}(\T)$.

The required $H^{2q}$-regularity of the Bloch waves $p_{n_j}(\cdot, k), j=1,\dots,N$ is satisfied if $V\in H^\mu(\T), \mu> 2q+d-2$, see Lemma 3 in \cite{DR18}. We obtain
$$
\begin{aligned}
\left\|\eps^{5/2-d}\pa_T\Atil_{n,l}\left(\frac{\cdot-l}{\eps},T\right)\right\|_{L^1(\B)}\leq &c\eps^{5/2}n^{-2(q+1)/d}\left(\sum_{j=1}^N \|\pa_T\Ahat_j(\cdot,T)\|_{L^1(\R^d)}\right. \\
& \left. +\sum_{\alpha,\beta,\gamma\in\{1,\dots,N\}} \|\pa_T\Ahat_\alpha(\cdot,T)\|_{L^1(\R^d)}\|\Ahat_\beta(\cdot,T)\|_{L^1(\R^d)}\|\Ahat_\gamma(\cdot,T)\|_{L^1(\R^d)}\right)
\end{aligned}
$$
and
\beq\label{E:dTAtil-est}
\begin{aligned}
\left\|\eps^{5/2-d}\left(\pa_T\Atil_{n,l}\left(\frac{\cdot-l}{\eps},T\right)\right)_{n\in\N}\right\|_{\cX(s)}\leq &c\eps^{5/2}\left(\sum_{j=1}^N \|\pa_T\Ahat_j(\cdot,T)\|_{L^1(\R^d)}\right. \\
& \hspace{-1.5cm}\left. +\sum_{\alpha,\beta,\gamma\in\{1,\dots,N\}} \|\pa_T\Ahat_\alpha(\cdot,T)\|_{L^1(\R^d)}\|\Ahat_\beta(\cdot,T)\|_{L^1(\R^d)}\|\Ahat_\gamma(\cdot,T)\|_{L^1(\R^d)}\right)
\end{aligned}
\eeq
provided $q> \tfrac{s}{2}+\tfrac{d}{4}-1$.

Next, we have
\beq\label{E:WU1-est}
\begin{aligned}
\left\|\eps \sum_{m\in \{-m_*,\dots, m_*\}}M^{(m)}(\cdot)\vv{U}^{\text{ext},1}(\cdot-l^{(m)},t)\right\|_{\cX(s)}\leq &c\eps \|W \tilde{u}^{\text{ext},1}(\cdot,\cdot,t)\|_{L^1(\B,H^s(\T))} \\
& \leq c\eps \|\tilde{u}^{\text{ext},1}(\cdot,\cdot,t)\|_{L^1(\B,H^s(\T))} \leq c\eps \|\vv{U}^{\text{ext},1}(\cdot,t)\|_{\cX(s)}
\end{aligned}
\eeq
For the nonlinear terms in $\hot$ we use the algebra property 
\beq\label{E:algeb_L1Hs}
\|\util *_\B \tilde{v}\|_{L^1(\B,H^s(\T))} \leq c \|\util\|_{L^1(\B,H^s(\T))}\|\tilde{v}\|_{L^1(\B,H^s(\T))}
\eeq
if $s>d/2$, see Lemma 2 in \cite{DR18}.
This yields, for instance
\beq\label{E:NLU1-est}
\left\|\vec{F}(\vv{U}^{\text{ext},1},\vv{U}^{\text{ext},0},\vv{U}^{\text{ext},0})\right\|_{\cX(s)} \leq c \|\vv{U}^{\text{ext},1}(\cdot,t)\|_{\cX(s)}\|\vv{U}^{\text{ext},0}(\cdot,t)\|_{\cX(s)}^2.
\eeq

To make inequalities \eqref{E:WU1-est} and \eqref{E:NLU1-est} useful, it remains to estimate $\|\vv{U}^{\text{ext},0}(\cdot,t)\|_{\cX(s)}$ and $\|\vv{U}^{\text{ext},1}(\cdot,t)\|_{\cX(s)}$. We have
\beq\label{E:U0-est}
\|\vv{U}^{\text{ext},0}(\cdot,t)\|_{\cX(s)} \leq c \eps^{1/2} \sum_{j=1}^N \|\Ahat_j(\cdot,T)\|_{L^1(\R^d)}
\eeq
and
\beq\label{E:U1-est}
\begin{aligned}
\|\vv{U}^{\text{ext},1}(\cdot,t)\|_{\cX(s)} &=\left(\sum_{n\in \N}n^{2s/d}\|U_n^{\text{ext},1}(\cdot,t)\|_{L^1(\B)}^2\right)^{1/2}\\
&\leq c\eps^{3/2} \left(\sum_{n\in \N}n^{2s/d-4(q+1)/d}\right)^{1/2}\left(\sum_{j=1}^N \|\Ahat_j(\cdot,T)\|_{L^1(\R^d)}\right. \\
& \left. +\sum_{\alpha,\beta,\gamma\in\{1,\dots,N\}} \|\Ahat_\alpha(\cdot,T)\|_{L^1(\R^d)}\|\Ahat_\beta(\cdot,T)\|_{L^1(\R^d)}\|\Ahat_\gamma(\cdot,T)\|_{L^1(\R^d)}\right)\\
&\leq c\eps^{3/2} \left(\sum_{j=1}^N \|\Ahat_j(\cdot,T)\|_{L^1(\R^d)}\right. \\
& \left. +\sum_{\alpha,\beta,\gamma\in\{1,\dots,N\}} \|\Ahat_\alpha(\cdot,T)\|_{L^1(\R^d)}\|\Ahat_\beta(\cdot,T)\|_{L^1(\R^d)}\|\Ahat_\gamma(\cdot,T)\|_{L^1(\R^d)}\right)\\
\end{aligned}
\eeq
if, again, $q> \tfrac{s}{2}+\tfrac{d}{4}-1$. In \eqref{E:U1-est} we used 
$\sigma\in H^{2q}(\T)$ and $V\in H^\mu(\T), \mu> 2q+d-2$, similarly to above, in order to get the decay in $n$.

In summary, collecting \eqref{E:Res-est1}, \eqref{E:dTAtil-est}, \eqref{E:WU1-est}, \eqref{E:NLU1-est}, \eqref{E:U0-est} and \eqref{E:U1-est} if $\sigma\in H^{2q}(\T), V\in H^\mu(\T), \mu> 2q+d-2$, $q\in \N\cap (\tfrac{s}{2}+\tfrac{d}{4}-1,\infty)$, and $\Ahat_j(\cdot,T)\in L^1_2(\R^2), \pa_T \Ahat_j(\cdot,T) \in L^1(\R^2)$ for all $T\in [0,T_0]$, then
\beq\label{E:Resr-est}
\|\vec{\text{Res}}(\cdot,t)\|_{\cX(s)}\leq c\eps^{5/2} \quad \text{for all } t \in [0,T_0\eps^{-1}],
\eeq
where $c>0$ depends polynomially on $\|\widehat{A}_j(\cdot,T)\|_{L^1_2(\R^d)}$ and  $\|\pa_T \Ahat_j(\cdot,T)\|_{L^1(\R^d)}, j=1,\dots,N$.

%----------------------------------------------------------------------
\subsection{Estimation of the Error}

Using the triangle inequality, we estimate the error $\|(u-u^\text{app})(\cdot,t)\|_{C^0_b(\R^d)}$ by the sum of $\|(\vv{U}^{\text{ext}}-\vv{U}^{\text{app}})(\cdot,t)\|_{\cX(s)}$ and $\|\vv{E}(\cdot,t)\|_{\cX(s)}$, where
$$\vv{E}:=\vv{U}-\vv{U}^{\text{ext}}.$$ 
We estimate $\|(\vv{U}^{\text{ext}}-\vv{U}^{\text{app}})(\cdot,t)\|_{\cX(s)}$ directly and $\|\vv{E}(\cdot,t)\|_{\cX(s)}$ via Gronwall's inequality

\blem\label{L:Uapp-Uext}
If $V\in H^{\mu}(\T)$ and $\sigma \in H^{2q}(\T)$ with some $\mu>2q+d-2$, $q\in \N\cap(\frac{s}{2}+\frac{d}{4},\infty)$, and $\Ahat_j(\cdot,T)\in L^1_\beta(\R^d)$, $j=1,\dots,N$, with some $\beta>2q+d$, then there is $c>0$ such that
$$
\begin{aligned}
\|(\Uvecapp-\Uvecext)(\cdot,t)\|_{\cX(s)}\leq &c\left(\eps^{3/2} \sum_{j=1}^N\|\Ahat_j(\cdot,T)\|_{L^1_\beta(\R^d)} \right. \\
&\left.+ \eps^{1/2+2d}\sum_{\alpha,\beta,\gamma\in\{1,\dots,N\}} \|\Ahat_\alpha(\cdot,T)\|_{L^1(\R^d)}\|\Ahat_\beta(\cdot,T)\|_{L^1(\R^d)}\|\Ahat_\gamma(\cdot,T)\|_{L^1(\R^d)}\right)
\end{aligned}
$$
for all $\eps>0$ small enough.
\elem
\bpf
Because $\Ahat_j(\cdot,T)\in L^2(\R^d)$, we have also $A_j(\cdot,T)\in L^2(\R^d)$ leading to $u^\text{app}(\cdot,t)\in L^2(\R^d)$. Therefore, it makes sense to evaluate $\cD\cT u^\text{app}(\cdot,t)$, producing
$$U^\text{app}_{n}(k,t)=\eps^{1/2-d}\sum_{j=1}^N\sum_{\eta\in\Z^d}\Ahat_j\left(\frac{k-k^{(j)}+\eta}{\eps},T\right)\langle p_{n_j}(\cdot,k^{(j)}-\eta),p_n(\cdot,k)\rangle_{L^2(\T)}, n \in\N.$$

It is
$$
\begin{aligned}
\|(\Uvecapp-\Uvecext)(\cdot,t)\|_{\cX(s)} \leq &c\sum_{j=1}^N\|(U^\text{app}_{n_j}-U^{\text{ext},0}_{n_j})(\cdot,t)\|_{L^1(\B)} + \|(U^\text{app}_{n}(\cdot,t))_{n\in \N\setminus\{n_1,\dots,n_N\}}\|_{\cX(s)}\\
&+ \|\vv{U}^{\text{ext},1}(\cdot,t)\|_{\cX(s)}.
\end{aligned}
$$

$\|\vv{U}^{\text{ext},1}(\cdot,t)\|_{\cX(s)}$ is estimated in \eqref{E:U1-est}. The remaining two contributions are estimated in a complete analogy to Lemma 10 in \cite{DR18}. We get
$$\|(U^\text{app}_{n_j}-U^{\text{ext},0}_{n_j})(\cdot,t)\|_{L^1(\B)} \leq c(\eps^{3/2}\|\Ahat_j(\cdot,T)\|_{L^1_1(\R^d)}+\eps^{1/2+\nu/2}\|\Ahat_j(\cdot,T)\|_{L^1_\nu(\R^d)})$$
for any $\nu>0$ using the Lipschitz continuity $|\langle p_{n_j}(\cdot,k^{(j)}),p_{n_j}(\cdot,k^{(j)})-p_{n_j}(\cdot,k)\rangle_{L^2(\T)}|\leq L |k-k^{(j)}|$ guaranteed by (A6).

Finally,
\beq\label{E:Unapp-est}
\|(U^\text{app}_{n}(\cdot,t))_{n\in \N\setminus\{n_1,\dots,n_N\}}\|_{\cX(s)} \leq c\left(\eps^{3/2} \sum_{j=1}^N\|\Ahat_j(\cdot,T)\|_{L^1_1(\R^d)} + \eps^{1/2+\beta}  \sum_{j=1}^N\|\Ahat_j(\cdot,T)\|_{L^1_\beta(\R^d)}\right)
\eeq
provided $\beta > 2q+d, q\in \N\cap (\tfrac{s}{2}+\tfrac{d}{4},\infty),$ $V\in H^{a}(\T), a>2q+d-2$. This estimate uses the Lipschitz continuity of $p_{n_j}(x,\cdot)$, the growth of the eigenvalues $\omega_n$, and the fact $\|\cL^q(k)p_{n_0}(\cdot,k-\eta)\|_{L^2(\T)} \leq c|\eta|^{2q}$. For details see Lemma 10 in \cite{DR18}.

The factor $\eps^{1/2+2d}$ in the Lemma comes from \eqref{E:Unapp-est} and  $\beta > 2q+d,q>s/2+d/4,s>d/2,$ which implies $\beta>2d$. 
\epf
Let us consolidate the regularity assumptions needed in the proof of Theorem \ref{T:justif}. The estimates are performed in $\cX(s), s>d/2$. The residual estimate \eqref{E:Resr-est} and Lemma \ref{L:Uapp-Uext} require $\hat{A}_j\in C([0,T_0],L^1_\beta(\R^d)\cap L^2(\R^d)), \pa_T\Ahat_j\in C([0,T_0],L^1(\R^d))$ for some $\beta > 2q+d$, $q\in \N\cap (\tfrac{s}{2}+\tfrac{d}{4},\infty)$. 
We can choose $q=\left\lceil\tfrac{d}{2}\right\rceil+1$ with $s\in (\tfrac{d}{2},2(\left\lceil\tfrac{d}{2}\right\rceil +1)-\tfrac{d}{2})$. Then 
the potentials need to satisfy $\sigma \in H^{2\left\lceil\tfrac{d}{2}\right\rceil+2}(\T)$ and $V\in H^{\mu}(\T)$ with $\mu>2\left\lceil \tfrac{d}{2}\right\rceil+d$. In summary, the regularity conditions are satisfied  if we require the existence of $\delta>0$ such that $\hat{A}_j\in C([0,T_0],L^1_{2\left\lceil\tfrac{d}{2}\right\rceil+d+2+\delta}(\R^d)\cap L^2(\R^d)),$ $\pa_T\Ahat_j\in C([0,T_0],L^1(\R^d))$, and $\sigma\in H^{2\left\lceil \tfrac{d}{2}\right\rceil+2}(\T), V\in H^{2\left\lceil \tfrac{d}{2}\right\rceil+d+\delta}(\T)$.

Finally, we estimate the error $\vec{E}$. The procedure is completely analogous to that in Sec. 4 of \cite{DH17}. We present it here for completeness.  The error satisfies
\beq\label{E:err_eq}
 \partial_t \Evec = -\ri \Omega(k) \Evec +\ri \vec{G}(\Uvecext,\Evec)
\eeq
with 
$$
\vec{G}(\Uvecext,\Evec)=\vec{\text{Res}}+\vec{F}(\Uvec,\Uvec,\Uvec) -\vec{F}(\Uvecext,\Uvecext,\Uvecext) -\varepsilon \sum_{m=-m_*}^{m_*} M^{(m)} \Evec(\cdot+l^{(m)},t).
$$
Due to the cubic structure of $\vec{F}$, and the isomorphism property of $\cD$ and the algebra property \eqref{E:algeb_L1Hs} 
we have for $s>d/2$ the existence of $c>0$ such that 
\beq\label{E:F_est}
\|\vec{F}(\Uvec,\Uvec,\Uvec) -\vec{F}(\Uvecext,\Uvecext,\Uvecext)\|_{\cX(s)} \leq c(\|\Uvecext\|^2_{\cX(s)}\|\Evec\|_{\cX(s)}+\|\Uvecext\|_{\cX(s)}\|\Evec\|^2_{\cX(s)}+\|\Evec\|^3_{\cX(s)}).
\eeq
We also have
\beq\label{E:M_est}
\left\| \varepsilon \sum_{m=-m_*}^{m_*} M^{(m)}(\cdot) \Evec(\cdot+l^{(m)},t)\right\|_{\cX(s)} \leq c \varepsilon \| \Evec(\cdot, t)\|_{\cX(s)}.
\eeq
Next, using \eqref{E:U0-est} and \eqref{E:U1-est}, we get
\beq\label{E:Uext_est}
\|\Uvecext(\cdot,t)\|_{\cX(s)}\leq c\eps^{1/2} \qquad \text{for all } t\in [0,\eps^{-1}T_0].
\eeq
Hence, there are $c_1,c_2,c_3>0$ such that
$$\|\vec{G}(\Uvecext,\Evec)(t)\|_{\cX(s)} \leq c_1\eps \|\Evec(t)\|_{\cX(s)}+c_2\eps^{1/2} \|\Evec(t)\|_{\cX(s)}^2+c_3\|\Evec(t)\|_{\cX(s)}^3+C_{\text{Res}}\eps^{5/2}$$
for all $t\in [0,\eps^{-1}T_0]$.

We rewrite equation \eqref{E:err_eq} as
$$
\Evec(t) = \Evec(0)+\int \limits_0^t S(t-\tau)\vec{G}(\Uvecext,\Evec)(\tau) ~d\tau,
$$
where $S(t)=e^{-\ri \Omega t}: \cX(s) \rightarrow \cX(s)$ is a strongly continuous unitary group  generated by $-\ri\Omega(k)$. 
Hence, we have 
$$
\begin{aligned}
\| \Evec(t) \|_{\cX(s)} \leq \| \Evec(0) \|_{\cX(s)}+\int \limits^t_0 c_1 \varepsilon \| \Evec(\tau) \|_{\cX(s)} + c_2 \varepsilon^{1/2} \| \Evec(\tau) \|^2_{\cX(s)}
+c_3 \| \Evec(\tau) \|^3_{\cX(s)} + C_{\text{Res}} \varepsilon^{5/2} ~d\tau.
\end{aligned}
$$
Because $\Uvec(0)=\Uvec^\text{app}(0)$, we get $\Evec(0)=\Uvec^\text{app}(0)-\Uvecext(0)$ and Lemma \ref{L:Uapp-Uext} provides $\|\Evec(0)\|_{\cX(s)}\leq C_0\eps^{3/2}$ for all $\eps\in (0,\eps_0)$ with some $\eps_0>0$.

Given an $M>C_0$ there exists $T_*>0$ such that $\|\Evec(t)\|_{\cX(s)}\leq M \eps^{3/2}$ for all $t\in [0,T_*]$. In the following, using Gronwall's inequality and a bootstrapping argument, we choose $\eps_0>0$ and $M>0$ such that if $\eps\in (0,\eps_0)$, then $\|\Evec(t)\|_{\cX(s)}\leq M \eps^{3/2}$ for all $t\in [0,\eps^{-1}T_0]$.

If $\|\Evec(t)\|_{\cX(s)} \leq M\eps^{3/2}$, then 
$$ 
\begin{aligned}
\| \Evec(t) \|_{\cX(s)} & \leq ~C_0 \eps^{3/2}+\int \limits^t_0 c_1 \varepsilon \| \Evec(\tau) \|_{\cX(s)}  ~d\tau + 
t\left( c_2 \varepsilon^{7/2} M^2 +c_3 \varepsilon^{9/2} M^3  +  C_{\text{Res}} \varepsilon^{5/2} \right)\\
& \leq ~ \varepsilon^{3/2} \left[ C_0 + t \eps\left( c_2 \varepsilon M^2 + c_3 \varepsilon^{2}  M^3 +  C_{\text{Res}} \right) \right] e^{c_1 \varepsilon t},
\end{aligned}
$$
where the second step follows by Gronwall's inequality.
For the desired estimate on $t\in [0,\eps^{-1}T_0]$, we redefine
$$M:=C_0+T_0(C_\text{Res}+1)e^{c_1T_0}$$
and choose $\eps_0$ so small that $ c_2 \varepsilon_0 M^2 + c_3 \varepsilon_0^{2}  M^3 \leq 1$. Then, clearly,
$$
\sup \limits_{t \in [0,\eps^{-1}T_0]} \| \Evec(t) \|_{\cX(s)} \leq M \varepsilon^{3/2}. 
$$
This completes the proof of Theorem \eqref{T:justif}.

%----------------------------------------------------------
\section*{Acknowledgments}
This research is supported by the \emph{German Research Foundation}, DFG grant No. DO1467/3-1.

\bibliographystyle{plain}
\bibliography{bib-CME-Rd}

\end{document}